\documentclass[smallextended,referee,envcountsect]{svjour3} 
\smartqed 
\usepackage{graphicx}
\usepackage[letterpaper, margin=1in]{geometry}

\usepackage{appendix}
\usepackage{amsfonts}
\usepackage{amsmath}
\usepackage{amssymb}
\usepackage{bbm}
\usepackage{bm}
\usepackage{mathrsfs}
\usepackage[inline, shortlabels]{enumitem}
\usepackage{verbatim}
\usepackage{setspace}
\usepackage{color}
\usepackage{pdfsync}
\usepackage{enumitem}
\usepackage{bbm}

\newtheorem{assumption}{Assumption}[section]
\newtheorem{re}{Remark}[section]

\newcommand{\eps}{\varepsilon}
\newcommand{\vp}{\varphi}

\newcommand{\N}{\mathbb{N}}

\newcommand{\R}{\mathbb{R}}

\renewcommand{\P}{\mathbb{P}}
\renewcommand{\L}{\mathbb{L}}

\newcommand{\bU}{\mathbb{U}}

\newcommand{\Uc}{\mathcal{U}}
\newcommand{\Di}{\mathbb{D}_{i}}
\newcommand{\D}{\mathbb{D}}
\newcommand{\DT}{\mathbb{D}_{T}}

\newcommand{\as}{\mbox{-a.s.}}

\def\be{\begin{eqnarray}}
\def\ee{\end{eqnarray}}
\def\beq{\begin{equation}}
\def\eeq{\end{equation}}

\journalname{JOTA}

\begin{document}

\title{Stochastic Perron for Stochastic Target Problems}

\author{Erhan Bayraktar  \and  Jiaqi Li}

\institute{Erhan Bayraktar,  Corresponding author  \at
             University of Michigan \\
              Ann Arbor, Michigan\\
              erhan@umich.edu  
           \and
           Jiaqi Li  \at
              University of Michigan \\
              Ann Arbor, Michigan\\
              lijiaqi@umich.edu
}

\date{Received: date / Accepted: date}

\maketitle

\begin{abstract}
In this paper, we adapt stochastic Perron's method to analyze stochastic target problems in a jump diffusion setup, where the controls are unbounded. Since classical control problems can be analyzed under the framework of stochastic target problems (with unbounded controls), we use our results to generalize the results of Bayraktar and S{\^{\i}}rbu (SIAM Journal on Control and Optimization, 2013) to problems with controlled jumps. 
\end{abstract}


\section{Introduction}\label{sec:intro}
Introduced by the seminal papers  \cite{Soner_Touzi_Superreplication}, \cite{DP_FOR_STP_AND_G_FLOW} and \cite{SONER_TOUZI_STG}, the stochastic target problem is a new type of optimal control problem. The aim is to drive a controlled diffusion to a given target at a pre-specified terminal time by choosing an appropriate admissible control. 
The above papers and their generalizations \cite{Bruno_jump_diffusion,Moreau} (to jump diffusions), \cite{Bouchard_Elie_Touzi_ControlledLoss} (to unbounded controls) provide a characterization of the associated value function as a viscosity solution to a non-linear  Hamilton-Jacobi-Bellman (HJB) equation using the geometric dynamic programming principle proved in \cite{DP_FOR_STP_AND_G_FLOW}.

In this paper, our goal is to provide an analysis of this problem using stochastic Perron's method.
This method  was introduced in \cite{Bayraktar_and_Sirbu_SP_LinearCase,Bayraktar_and_Sirbu_SP_DynkinGames,Bayraktar_and_Sirbu_SP_HJBEqn} for classical control problems. This method is a verification approach (without requiring smoothness) in that it does not use the dynamic programming principle to show that the value function is a viscosity solution. The idea is to build two classes of functions that envelope the value function and that are stable enough under minimization and maximization, respectively. This construction helps us demonstrate that the supremum over the first class is a lower semi-continuous viscosity super-solution and the infimum over the second class (the functions larger than the value function) is an upper semi-continuous  viscosity sub-solution. Assuming that a comparison principle holds,  we show that the infimum over the second class and the supremum over the first class (which sandwich the value function) are equal, and hence, the value function is the unique viscosity solution. Since we only work with the envelopes, not the value function itself,  we never use the dynamic programming principle (and hence the measurable selection theorem). In fact, the dynamic programming principle is a corollary of our result. As pointed out by \cite{MR3488161} and the references therein, the rigorous proof of the dynamic programming principle for controlled diffusion processes is difficult and contains subtle technical issues. Our result can be seen as an elementary alternative based only on It\^{o}'s Lemma and the comparison principle, which also has to be proved to identify the value function as the unique viscosity solution of the Hamilton-Jacobi-Bellman partial differential equation (PDE).

We choose to work with the most general stochastic target setup from
 \cite{Moreau}. Our controls are unbounded and the controlled processes are jump diffusions.  The main reason for using unbounded controls is that it allows us to use the embedding result of \cite{Bouchard_Equivalence}, which converts an ordinary control problem into a stochastic target problem with unbounded admissible controls. Using this result, we generalize \cite{Bayraktar_and_Sirbu_SP_HJBEqn} to the setting of controlled jumps. 
 
In contrast to \cite{Bayraktar_and_Sirbu_SP_HJBEqn}, we analyze stochastic target problems in this paper.  The main contribution is the construction of the sets of stochastic semi-solutions, which are appropriate for stochastic target problems. This makes the proofs of the viscosity properties of the value function different. We also generalize our earlier result in \cite{BayraktarLi} in the sense that we consider unbounded controls and controlled jumps. The presence of the jumps and the unbounded control set brings new technical difficulties: in contrast to \cite{BayraktarLi}, the relaxed semi-limits are introduced for the PDE characterization, which have a nontrivial impact on the formulation of the associated PDEs and the derivation of viscosity properties of the value function using stochastic Perron's method, especially at the boundary. Of particular importance is the relaxation with respect to the test function, which appears because we consider jumps.

The rest of the paper is organized as follows. The setup of the problem, the related HJB equation and the definitions of the stochastic semi-solutions are introduced in Section \ref{sec:prob}. In Sections \ref{sec: interior} and \ref{sec: bdd cond}, we prove the viscosity properties in the parabolic interior and at the boundary, respectively. In Section \ref{sec: Verification}, we use the comparison principle to close the gap between the viscosity super-solution and sub-solution and demonstrate the uniqueness of the viscosity solution to the associated HJB equation. In Section \ref{sec:equivalence}, we see how an optimal control problem can be converted into a stochastic target problem. Some technical results are delegated to the Appendix. Our main results are Theorems \ref{thm: main theorem_interior}, \ref{thm: bd_viscosity_property}, \ref{thm: unique viscosity solution} and \ref{thm: optimal control}.

\section{The Setup}\label{sec:prob}
To introduce the stochastic target problem in \eqref{eq: value_function}, we need to introduce some notation and make appropriate 
assumptions. Throughout this paper, the superscript $^{\top}$ stands for transposition, $|\cdot|$ for the Euclidean norm of a vector in $\R^n$ and $\|\cdot\|$ for the Frobenius norm of a matrix.  For a subset of $\mathcal{O}$ of $\R^n$, we denote by Int$(\mathcal{O})$ its interior. We also denote the open ball of radius $r>0$ centered at $x\in \R^n$ by $B_r(x)$ and the set of $n\times n$ matrices by $\mathbb{M}^n$.  Inequalities and inclusion between random variables and random sets, respectively, are in the almost sure sense unless otherwise stated.  

Given a complete probability space $(\Omega, \mathcal{F}, \mathbb{P})$, let $\{\lambda_i(\cdot,de)\}_{i=1}^I$ be a collection of independent integer-valued $E$-marked right-continuous point processes defined on this space. Here, $E$ is a Borel subset of $\mathbb{R}$ equipped with the Borel sigma field $\mathcal{E}$. Let $\lambda=(\lambda_1,\lambda_2,\cdots,\lambda_I)^{\top}$and $W = \{W_s\}_{0\leq s\leq T}$ be a $d$-dimensional Brownian motion defined on the same probability space such that $W$ and $\lambda$ are independent. Given $t\in[0,T]$, let $\mathbb{F}^t=\{\mathcal{F}^t_s, t\leq s\leq T\}$ be $\mathbb{P}$-completed filtration generated by $W_{\cdot}-W_t$ and $\lambda([0,\cdot],de)-\lambda([0,t],de)$.  Set $\mathcal{F}^t_s = \mathcal{F}^t_t $ for $0\leq s < t$. We will use $\mathcal{T}_t$ to denote the set of $\mathbb{F}^t$-stopping times valued in $[t,T]$. Given $\tau\in\mathcal{T}_t$, the set of $\mathbb{F}^t$-stopping times valued in $[\tau, T]$ will be denoted by $\mathcal{T}_{\tau}$.  
\begin{assumption}\label{assump: lambda intensity kernel}
$\lambda$ satisfies the following:
\begin{enumerate}
\item $\lambda(ds,de)$ has intensity kernel $m(de)ds$ such that $m_{i}$ is a Borel measure on $(E,\mathcal{E})$ for any $i=1,\cdots, I$ and $\hat{m}(E)<\infty$, where $m=(m_1,\cdots,m_I)^{\top}$ and $\hat{m}=\sum_{i=1}^{I} m_{i}$.
 \item $E=\text{supp}(m_i)$ for all $i=1,2,\cdots,I$. Here, $\text{supp}(m_{i}):=\{e\in E:  e\in N_{e}\in T_{E} \implies m_{i}(N_{e})>0\},$
 where $T_{E}$ is the topology on $E$ induced by the Euclidean topology. 
 \item There exists a constant $C>0$ such that 
 $$\mathbb{P}\left(\left\{\hat{\lambda}(\{s\}, E)\leq C \;\;\text{for all}\;\; s\in[0,T]\right\}\right)=1,\;\;\text{where}\;\;\hat{\lambda}=\sum_{i=1}^{I}\lambda_{i}.$$ 
\end{enumerate}
\end{assumption}
The above assumption implies that there are a finite number of jumps during any finite time interval. Let $\tilde{\lambda}(ds,de):=\lambda(ds,de)-m(de)ds$ be the associated compensated random measure.  

Let $\mathcal{U}^t_1$ be the collection of all the $\mathbb{F}^{t}$-predictable processes in $\mathbb{L}^2(\Omega\times[0,T], \mathcal{F}\otimes\mathcal{B}[0,T], \mathbb{P}\otimes \lambda_{L}; U_{1})$, where $\lambda_{L}$ is the Lebesgue measure on $\R$ and $U_{1}\subset \R^{q}$ for some $q\in\N$. Define $\mathcal{U}_2^{t}$ to be the collection of all the maps $\nu_{2}:\Omega\times[0,T]\times E\rightarrow \R^{n}$ which are $\mathcal{P}^{t}\otimes\mathcal{E}$ measurable such that 
\begin{equation*}
\|\nu_{2}\|_{\mathcal{U}_2^{t}}:=\left(\mathbb{E}\left[\int_{t}^{T}\int_E|\nu_{2}(s,e)|^2 \hat{m}(de)ds\right]\right)^{\frac{1}{2}}<\infty,
\end{equation*}
 where $\mathcal{P}^{t}$ is the $\mathbb{F}^t$-predictable sigma-algebra on $\Omega\times[0,T]$. $\nu=(\nu_{1},\nu_{2})\in\mathcal{U}_{0}^{t} := \mathcal{U}_{1}^{t}\times\mathcal{U}^{t}_{2}$ takes value in the set $U:=U_1\times \mathbb{L}^2(E,\mathcal{E}, \hat{m};\R^{n})$.  
Let $\D=[0,T]\times \R^d$, $\Di=[0,T[\;\times\;\R^d \text{ and } \DT=\{T\}\times \R^d.$
Given $z = (x, y)\in\R^d\times\R $, $t \in [0, T]$ and $\nu\in\mathcal{U}_{0}^{t}$, we consider the stochastic differential equations (SDEs)
\begin{equation}\label{eq: SDEs}
\begin{array}{l}
dX(s)=\mu_{X}(s,X(s),\nu(s))ds+\sigma_{X}(s,X(s),\nu(s))dW_s+\int_{E} \beta(s,X(s-),\nu_1(s),\nu_2(s,e), e)\lambda(ds,de), \vspace{0.07in}\\
dY(s)=\mu_{Y}(s,Z(s),\nu(s))ds+\sigma_{Y}^{\top}(s,Z(s),\nu(s)) dW_s+ \int_{E} b^{\top}(s,Z(s-),\nu_1(s),\nu_2(s,e), e)\lambda(ds,de), 
\end{array}
\end{equation}
with $(X(t), Y(t))=(x,y)$. Here, $Z=(X,Y)$. 
In \eqref{eq: SDEs}, 
 \begin{equation*}
\begin{array}{c}
\mu_X: \D\times U \rightarrow \R^d,\;\;\sigma_X: \D\times U\rightarrow \R^{d\times d},\;\; \beta: \D\times U_1\times\R^n\times E \rightarrow \R^{d\times I}, \\
\mu_Y: \D\times\R\times U \rightarrow \mathbb{R},\;\;\sigma_Y: \D\times\R\times U\rightarrow \R^{d},\;\; b: \D\times\R\times U_1\times\R^n\times E \rightarrow \R^{I}. 
\end{array}
\end{equation*}
Besides the measurability and the integrability conditions for $\mathcal{U}^{t}_{0}$, we impose another condition on the admissible control set. Let $\mathcal{U}^{t}$ be the admissible control set, which consists of all $\nu \in \mathcal{U}^{t}_{0}$ such that for any compact set $C\subset \R^{d}\times \R$, there exists a constant $K_{C, \nu}>0$ such that 
\begin{equation}\label{eq:admissibility}
\left|\int_{E} b^{\top}(\tau,x, y, \nu_{1}(\tau), \nu_{2}(\tau, e), e ) \lambda(\{\tau\}, e) \right| \leq K_{C, \nu}\;\;\text{for all}\;\; (x,y)\in C\;\;\text{and}\;\;\tau\in\mathcal{T}_{t}.
\end{equation}

\begin{assumption}\label{assump: regu_on_coeff}
Let $z=(x,y)$ and $u=(u_1,u_2)\in U = U_1\times\mathbb{L}^2(E,\mathcal{E},\hat{m};\R^{n})$. We use the notation $\|u\|_{U}:=|u_1|+\|u_2\|_{\hat{m}}$ and  $u(e):=(u_1,u_2(e))$ for the rest of the paper. 
\begin{enumerate}
\item  $\mu_X, \sigma_X$, $\mu_Y$ and $\sigma_Y$ are all continuous; 
\item  $\mu_X, \sigma_X$, $\mu_Y$, $\sigma_Y$ are Lipschitz in $z$ and locally Lipschitz in other variables. In addition,
\begin{equation*}
|\mu_X(t,x,u)|+|\sigma_X(t,x,u)|\leq L(1+|x|+\|u\|_{U}), \;\;|\mu_Y(t,x,y,u)|+|\sigma_Y(t,x,y,u)|\leq L(1+|y|+\|u\|_{U}).
\end{equation*} 
\item $b$ and $\beta$ are Lipschitz and grow linearly in all variables except $e$, but uniformly in $e$.
\end{enumerate}
\end{assumption}
\begin{remark} \label{remark: facts from assumption 1}
Assumptions \ref{assump: lambda intensity kernel} and \ref{assump: regu_on_coeff} guarantee that there exists a unique strong solution $(X_{t,x}^{\nu}, Y_{t,x,y}^{\nu})$ to \eqref{eq: SDEs} for any $\nu\in\Uc^{t}$. Moreover, the processes $(X_{t,x}^{\nu}$, $ Y_{t,x,y}^{\nu})$ are c\`adl\`ag.
\end{remark}
\begin{remark}
Under Assumptions \ref{assump: lambda intensity kernel} and \ref{assump: regu_on_coeff}, $\mathcal{U}^{t}$ contains all the bounded processes in $\mathcal{U}^{t}_{0}$.\footnote{The bound may depend on the process.}
\end{remark}
We now define the value function of the stochastic target problem. Let $g: \R^d\rightarrow \R$ be a measurable function with polynomial growth. The value function of the target problem is defined by 
\begin{equation}\label{eq: value_function}
u(t,x):=\inf\left\{y:\; \exists \nu \in \Uc^t \text{ s.t. } Y_{t,x,y}^{\nu}(T)\geq g(X_{t,x}^{\nu}(T)) \; \mathbb{P}-\text{a.s.}\right\}.
\end{equation}
\subsection{\bf The Hamilton-Jacobi-Bellman Equation}\label{subsection: HJB}
Denote $b=(b_1,b_2,\cdots,b_I)^{\top}$ and $\beta=(\beta_1,\beta_2,\cdots,\beta_I)$. For a given $\vp \in C(\D)$, we define the relaxed semi-limits 
\begin{equation}\label{eq: HJB operators}
H^{*}(\Theta, \vp):=\limsup_{\begin{subarray}{c}\eps\searrow 0,\; \Theta^{'}\rightarrow\Theta \\ \eta\searrow 0,\;  \psi\overset{\text{u.c.}}{\longrightarrow} \vp\end{subarray}} H_{\eps,\eta} (\Theta^{'}, \psi) \;\; \text{and} \;\; H_{*}(\Theta, \vp):=\liminf_{\begin{subarray}{c}\eps\searrow 0,\; \Theta^{'}\rightarrow\Theta \\  \eta\searrow 0,\; \psi \overset{\text{u.c.}}{\longrightarrow} \vp\end{subarray}} H_{\eps, \eta} (\Theta^{'}, \psi). \footnote{The convergence $\psi \overset{\text{u.c.}}{\longrightarrow} \vp$ is understood in the sense that $\psi$ converges uniformly on compact subsets to $\vp$.}
\end{equation}
Here, for $\Theta = (t,x,y,p,A)\in\D\times\R\times\R^d\times\mathbb{M}^d$, $\vp\in C(\D)$, $\eps\geq 0$ and $\eta\in[-1,1]$, 
\begin{equation*}
\begin{gathered}
H_{\eps, \eta}(\Theta, \vp):=\sup_{u\in\mathcal{N}_{\eps, \eta}(t,x,y,p,\vp)} \mathbf{F}^u(\Theta), \;\text{where},\\
\begin{array}{c}
\mathbf{F}^u(\Theta):=\mu_{Y}(t,x,y,u) -\mu_{X}^{\top}(t,x,u) p -\frac{1}{2}\text{Tr}[\sigma_X\sigma_X^{\top}(t,x,u)A],\;\;N^{u}(t,x,y,p):=\sigma_Y(t,x,y,u)-\sigma_X^{\top}(t,x,u)p, \\
\Delta^{u,e}(t,x,y,\vp):=\min_{1\leq i \leq I} \{b_i(t,x,y,u(e),e)-\vp(t,x+\beta_i(t,x,u(e),e))+\vp(t,x) \}, \\
\mathcal{N}_{\eps,\eta}(t,x,y,p, \vp):=\{u\in U: |N^{u}(t,x,y,p)|\leq\eps\text{ and }
\Delta^{u,e}(t,x,y,\vp)\geq\eta\;\text{for}\; \hat{m}-\text{a.s.}\; e\in E \}.
\end{array}
\end{gathered}
\end{equation*}
For our later use, we also define the following:
\begin{equation*}
\begin{array}{c}
J^{u,e}_i(t,x,y,\vp):= b_i(t,x,y,u(e),e)-\vp(t,x+\beta_i(t,x,u(e),e))+\vp(t,x), \nonumber \\
\overline{J}^{u,e}(t,x,y,\vp):= (J^{u,e}_1(t,x,y,\vp), \cdots, J^{u,e}_I(t,x,y,\vp))^{\top},\;  J^{u}(t,x,y, \vp):=\inf_{e\in E}\min_{1\leq i\leq I} J_i^{u,e}(t,x,y,\vp), \nonumber \\
 \mathscr{L}^{u}\vp(t,x):=\vp_t(t,x)+ \mu_{X}^{\top}(t,x,u)D\vp(t,x)+\frac{1}{2}\text{Tr}[\sigma_X\sigma_X^{\top}(t,x,u)D^2\vp(t,x)]. 
 \end{array}
  \end{equation*}
\begin{remark}
For simplicity, we denote $H^*(t,x,\vp(t,x),D\vp(t,x), D^2\vp(t,x),\vp)$ by $H^*\vp(t,x)$ for $\vp\in C^{1,2}(\D)$. For $\vp\in C^2(\R^d)$, we denote $H^*(T,x,\vp(x),D\vp(x), D^2\vp(x),\vp)$ by $H^*\vp(x)$. We will use similar notation for $H_*$ and other operators in later sections.
\end{remark}
Later, we will produce a viscosity super-solution and sub-solution, respectively, to
 \begin{gather}\label{eq: super_HJB equation_interior}
-\partial_t\vp(t,x)+ H^*\vp(t,x)\geq 0 \;\;\text{in}\;\;\Di\;\;\text{and} \\
\label{eq: sub_HJB equation_interior}
 -\partial_t\vp(t,x)+H_*\vp(t,x)\leq 0 \;\;\text{in}\;\;\Di.
 \end{gather} 

\subsection{\bf Stochastic Semi-Solutions}
Before we introduce the definitions of the stochastic semi-solutions, we define the concatenation of the admissible controls. 
\begin{definition}[Concatenation]
Let $\nu_{1}, \nu_{2} \in \Uc^t$, $\tau \in \mathcal{T}_t$. The concatenation of $\nu_1$ and $ \nu_2 $ at  $\tau$ is defined as $
\nu_1\otimes_{\tau} \nu_2 := \nu_1 \mathbbm{1}_{[0,\tau[}+ \nu_2 \mathbbm{1}_{[\tau,T]} \in\Uc^{t}$.\footnote{This can be easily checked.}
\end{definition}  
\begin{definition} [Stochastic Super-solutions] \label{def: Stochasticsuper-solution}
A continuous function $w: \D \rightarrow \mathbb{R}$ is called a stochastic super-solution if
\begin{enumerate}
\item $w(T, x)\geq g(x)$ and for some $C>0$ and $n\in\N$,\footnote{$C$ and $N$ may depend on $w$ and $T$. This also applies to Definition \ref{def: Stochasticsub-solution}} $|w(t,x)|\leq C(1+|x|^{n})$ for all $(t,x)\in \D$.
\item Given $(t,x,y)\in \D\times\mathbb{R}$, for any $\tau\in\mathcal{T}_t$ and $\nu\in \Uc^t$, there exists $\tilde {\nu}\in \Uc^t$ such that  
$Y(\rho )\geq w(\rho, X(\rho )) \ \ \mathbb{P}-\text{a.s.}$ on $\{ Y(\tau)\geq w(\tau, X(\tau)) \}$ for all $\rho \in \mathcal{T}_{\tau}$,
where $X:= X_{t,x}^{\nu\otimes_{\tau}\tilde{\nu}}$ and $Y:=Y_{t,x,y}^{\nu\otimes_{\tau}\tilde{\nu}}$.
\end{enumerate}
\end{definition}
\begin{definition} [Stochastic Sub-solutions] \label{def: Stochasticsub-solution}
A continuous function $w: \D \rightarrow \mathbb{R}$ is called a stochastic sub-solution if
\begin{enumerate}
\item $w(T, x)\leq g(x)$ and for some $C>0$ and $n\in\N$, $|w(t,x)|\leq C(1+|x|^{n})$ for all $(t,x)\in \D$.
\item Given $(t,x,y)\in \D\times\mathbb{R}$, for any $\tau\in\mathcal{T}_t$ and $\nu\in \Uc^t$, we have
$\mathbb{P}(Y(\rho )< w(\rho, X(\rho ))|B)>0$
for all $\rho\in \mathcal{T}_{\tau}$ and $B\subset \{Y(\tau)<w(\tau,X(\tau))\}$ satisfying $B\in\mathcal{F}_\tau^t$ and $\mathbb{P}(B)>0$. Here, we use the notation $X:= X_{t,x}^{\nu}$ and $Y:=Y_{t,x,y}^{\nu}$.
\end{enumerate}
\end{definition}
Denote the sets of stochastic super-solutions and sub-solutions by $\bU^+$ and $\bU^-$, respectively.
\begin{assumption}\label{assump:semisolution_not_empty}
$\bU^+$ and $\bU^-$ are not empty. 
\end{assumption}
\begin{remark}
Let $u^+:=\inf_{w \in \mathbb{U}^{+}} w$. For any stochastic super-solution $w$, choose $\tau=t$ and $\rho=T$. Then there exists $\tilde{\nu}\in \Uc^t$ such that
$Y_{t,x,y}^{\tilde{\nu}}(T)\geq w\left(T, X_{t,x}^{\tilde{\nu}}(T)\right)\geq g\left(X_{t,x}^{\tilde{\nu}}(T)\right) \mathbb{P}-\text{a.s.}
 \text{ if }  y\geq w(t,x).$
Hence, $y\geq w(t,x)$ implies that $y\geq u(t,x)$ from \eqref{eq: value_function}. This means that $w\geq u$ and $u^+\geq u$. By the definition of $\bU^{+}$, we know that $u^+(T,x)\geq g(x)$ for all $x\in\R^d$.
\end{remark}
\begin{remark}
Let $u^-:=\sup_{w \in \mathbb{U}^{-}} w$. For any stochastic sub-solution $w$, if $y<w(t,x)$,  by choosing $\tau=t$ and $\rho=T$, we get that for any $\nu\in \Uc^t$,
$\mathbb{P}\left(Y_{t,x,y}^{\nu}(T)<  g(X_{t,x}^{\nu}(T))\right)\geq \mathbb{P}\left(Y_{t,x,y}^{\nu}(T)<  w(T,X_{t,x}^{\nu}(T))\right)>0.
$
Therefore, from \eqref{eq: value_function}, $y<w(t,x)$ implies that $y\leq u(t,x)$. This means that $w\leq u$ and $u^-\leq u$.  By the definition of $\bU^{-}$, it holds that $u^{-}(T,x)\leq g(x)$ for all $x\in\R^d$.
\end{remark}
In short,
\begin{equation}\label{eq:intfvmavp}
u^- = \sup _{w\in \mathbb{U}^-} w\leq  u \leq \inf _{w\in \mathbb{U}^+}w= u^+.
\end{equation}
We will provide sufficient conditions which guarantee Assumption \ref{assump:semisolution_not_empty} in the Appendix A. As in \cite{Bruno_jump_diffusion} and \cite{Moreau}, the proof of the sub-solution property requires a regularity assumption on the set-valued map $\mathcal{N}_{0,\eta}(\cdot, \psi)$.
\begin{assumption}\label{assump: regularity}
For $\psi\in C(\D)$, $\eta>0$, let $B$ be a subset of $\D\times\R\times\R^d$ such that $\mathcal{N}_{0,\eta}(\cdot, \psi)\neq\emptyset$ on $B$. Then for every $\eps>0$, $(t_0,x_0,y_0,p_0)\in Int(B)$ and $u_0\in\mathcal{N}_{0, \eta}(t_0,x_0,y_0,p_0,\psi) $, there exists an open neighborhood $B'$ of $(t_0,x_0,y_0,p_0)$ and a locally Lipschitz continuous map $\hat{\nu}$ defined on $B'$ such that $\|\hat{\nu}(t_0,x_0,y_0,p_0)-u_0\|_{U}\leq \eps$ and $\hat{\nu}(t,x,y,p)\in \mathcal{N}_{0, \eta}(t,x,y,p, \psi)$.
\end{assumption}  

\section{Viscosity Property in $\Di$}\label{sec: interior}
In this section, we state and prove the theorem which characterizes $u^+$ (resp. $u^-$) as a viscosity sub-solution (resp. super-solution) of \eqref{eq: sub_HJB equation_interior} (resp. \eqref{eq: super_HJB equation_interior}). The boundary conditions will be discussed in Theorem \ref{thm: bd_viscosity_property}.
\begin{lemma}
$\bU^{+}$ and $\bU^{-}$ are closed under pairwise minimization and maximization, respectively. That is, \\
\begin{enumerate*}[series = tobecont, itemjoin = \;\;]
\item if $w_1, w_2\in \mathbb{U}^+$, then $w_1\wedge w_2\in \mathbb{U}^+$;\; \item if $w_1,w_2\in \mathbb{U}^-$, then $w_1\vee w_2\in \mathbb{U}^-$.
\end{enumerate*}
\end{lemma}
\begin{lemma}\label{lem: monotone seq approaches v+ or v_-}
There exists a non-increasing sequence $\{w_{n}\}_{n=1}^{\infty}\subset\bU^{+}$ such that $w_n\searrow u^+$ and a non-increasing sequence $\{v_{n}\}_{n=1}^{\infty}\subset\bU^{-}$ such that $v_n\nearrow u^-$.
\end{lemma}
\begin{theorem} \label{thm: main theorem_interior}
Under Assumptions \ref{assump: lambda intensity kernel}-\ref{assump: regularity}, 
$u^+$ is an upper semi-continuous (USC) viscosity sub-solution of \eqref{eq: sub_HJB equation_interior}. On the other hand, under Assumptions \ref{assump: lambda intensity kernel}-\ref{assump:semisolution_not_empty}, $u^-$ is a lower semi-continuous (LSC) viscosity super-solution of \eqref{eq: super_HJB equation_interior}.
\end{theorem}
{\it Proof } See Appendix B. \qed

\section{Boundary Conditions}\label{sec: bdd cond}
\noindent In this section, we discuss the boundary conditions at $T$. From the definition of the value function $u$, it holds that $u(T,x)=g(x)$ for all $x\in\R^d$. However, $u^+$ and $u^-$ may not satisfy this boundary condition. Define
\begin{equation*}
\mathbf{N}(t,x,y,p, \psi):=\{(r,s)\in\R^d\times\R: \exists u \in U, \;\text{s.t.}\; r=N^{u}(t,x,y,p)\;\text{and}\;s\leq \Delta^{u,e}(t,x,y,\psi)\; \hat{m}-\text{a.s.}\}
\end{equation*}
and $
\delta:=\text{dist}(0, \mathbf{N}^c)-\text{dist}(0, \mathbf{N}),$
where dist denotes the Euclidean distance. It holds that
\begin{equation}\label{eq: delta>0 equi int}
0\in\text{int}(N(t,x,y,p,\psi))\;\;\text{iff}\;\;\delta(t,x,y,p,\psi)>0. 
\end{equation}
 The upper (resp. lower) semi-continuous envelope of $\delta$ is denoted by $\delta^* \;(\text{resp}.\; \;\delta_*)$. 
Let
\begin{equation*}\label{eq: boundary_u-_def}
u^+(T-,x)=\limsup _{(t<T, x')\rightarrow (T,x)} u^-(t,x'), \;\;u^-(T-,x)=\liminf _{(t<T, x')\rightarrow (T,x)} u^-(t,x').
\end{equation*}
The following theorem is an adaptation of the results in \cite{DP_FOR_STP_AND_G_FLOW,SONER_TOUZI_STG,Bruno_jump_diffusion,Bouchard_Equivalence}.
\begin{theorem}\label{thm: bd_viscosity_property}
Under Assumptions \ref{assump: lambda intensity kernel}-\ref{assump: regularity}, if $g$ is USC, then $u^{+}(T-,\cdot)$ is a USC viscosity sub-solution of 
$\min\{\vp(x)-g(x), \delta_{*}\vp(x)\}\leq 0\text{ on } \R^d.$
On the other hand, under Assumptions \ref{assump: lambda intensity kernel}-\ref{assump:semisolution_not_empty}, if $g$ is LSC, $u^{-}(T-,\cdot)$ is an LSC viscosity super-solution of 
$
\min\{(\vp(x)-g(x))\mathbbm{1}_{\{H^{*}\vp(x)<\infty\}}, \delta^{*}\vp(x) \}\geq 0\text{ on }\R^d.
$
\end{theorem}
{\it Proof } {\bf Step 1 (The sub-solution property on $\DT$).}
For the sake of contradiction, we assume that for some $x_0\in\R^d$ and $\vp\in C^{2}(\R^d)$ satisfying 
$0=u^+(T-,x_0)-\vp(x_0)= \max_{x\in \R^d}(u^{+}(T-,x)-\vp(x)),$
it holds that $\vp(x_0)-g(x_0)>2\eta$ and $\delta_{*}\vp(x_0)>2\eta \text{ for some } \eta>0.$
Let $\{w_{k}\}_{k=1}^{\infty}$ be a sequence in $\mathbb{U}^+$ such that $w_k\searrow u^+$. Set $\tilde{\vp}(t,x)=\vp(x)+\iota|x-x_0|^{n_{0}}+\iota\sqrt{T-t}$ for $\iota>0$, where $\iota$ will be fixed later and $n_{0}$ satisfies
\begin{equation*}
\min_{0\leq t\leq T}(\tilde{\vp}(t,x)-w_{1}(t,x))\rightarrow\infty\;\;\text{as}\;\;|x|\rightarrow\infty\;\;\text{for any}\;\; \iota>0. 
\end{equation*}
 By the lower semi-continuity of $\delta_*$ and the upper semi-continuity of $g$, we can find  $\iota>0$ and $\eps>0$ such that
\begin{gather}
\label{eq: vp great than g}
\tilde\vp(t,x)-g(x)>\eta \;\;\text{and} \\
\label{eq:delta great than 0 in the neighborhood}
\delta_{*}(t,x,y,D\tilde{\vp}(t,x),\tilde{\vp})\geq \eta \;\;
\text{for}\;\;(t,x)\in[T-\eps, T]\times \text{cl}(B_{\eps}(x_0)) \;\;\text{and}\;\; |y-\tilde{\vp}(t,x)|\leq \eps.
\end{gather}
 By Assumption \ref{assump: regularity}, the fact that $\delta\geq \delta_*$, \eqref{eq: delta>0 equi int} and \eqref{eq:delta great than 0 in the neighborhood}, we can find  a locally Lipschitz map $\hat{\nu}$ such that
\begin{equation}\label{eq: Lipschitz map exists_u+_bdd}
\begin{array}{c}
\hat{\nu}(t,x,y,D\tilde{\vp}(t,x))\in \mathcal{N}_{0,\eta}(t,x,y,\tilde{\vp}(t,x),\tilde{\vp}) \\
\text{for all}\;\; (t,x,y)\in \D\times R \;\;\text{s.t.}\;\; (t,x)\in[T-\eps, T]\times \text{cl}(B_{\eps}(x_0)) \;\;\text{and}\;\; |y-\tilde{\vp}(t,x)|\leq \eps.
\end{array}
\end{equation} 
In \eqref{eq: Lipschitz map exists_u+_bdd}, we may need to choose smaller values of $\eps,\iota$ and $\eta$. Fix $\iota$. Since $\partial_t\tilde{\vp}(t,x)\rightarrow -\infty$ as $t\rightarrow T$, by the continuity of $\mu_Y, \mu_X,\sigma_{X}$ and $\nu$, 
\begin{equation}\label{eq: drift great than 0 boundary u+}
\begin{array}{c}
\mu_Y(t,x,y,\hat\nu(t,x,y,D\tilde{\vp}(t,x)))-\mathcal{L}^{\hat\nu(t,x,y,D\tilde{\vp}(t,x))}\tilde{\vp}(t,x)\geq \eta,  \\
\text{for all}\;\; (t,x,y)\in \D\times R \;\;\text{s.t.}\;\; (t,x)\in[T-\eps, T]\times \text{cl}(B_{\eps}(x_0)) \;\;\text{and}\;\; |y-\tilde{\vp}(t,x)|\leq \eps.
\end{array}
\end{equation}
Here we may need to shrink $\eps>0$ again. Since $u^+$ is USC and $\tilde\vp(T,x_0)=u^+(T-,x_0)$, there exists $\alpha>0$ such that $\tilde{\vp}>u^+-2\alpha$ on  $[T-\eps, T[\;\times\;\text{cl}(B_{\eps/2}(x_0))$  after possibly shrinking $\eps$ another time. Since $w_k\searrow u^+$, there exists $n_0\in\N$ such that 
\begin{equation}\label{eq:comp1_boundary_u+}
\tilde\vp>w_{n_0}-\alpha\;\;\text{on} \;\;[T-\eps,T[\;\times\;\text{cl}(B_{\eps/2}(x_0)).
\end{equation}
Since $\min_{0\leq t \leq T}(\tilde\vp(t,x)-w_{1}(t,x))\rightarrow \infty$ as $|x|\rightarrow\infty$, we can find $R_0>\eps$ such that
\begin{equation}\label{eq:comp2_boundary_u+}
\tilde{\vp}>w_{n_0}+\eps \text{ on } \mathbb{O}:=[T-\eps, T]\times(\R^d\setminus\text{cl}(B_{R_0}(x_0))). 
\end{equation}
Notice that $\tilde\vp(T,\cdot)-u^{+}(T-,\cdot)$ is strictly positive on the compact set $\mathbb{T}^{*}:=\text{cl}(B_{R_{0}}(x_0))-B_{\eps/2}(x_0)$. Hence, by the upper semi-continuity of $u^{+}(T-,\cdot)$, there exists $\zeta>0$ such that
\begin{equation}\label{eq:vp great than u+ on DT}
\tilde\vp(T,\cdot)>u^{+}(T-,\cdot)+4\zeta\;\;\text{on}\;\; \mathbb{T}^{*}.
\end{equation}
From \eqref{eq:vp great than u+ on DT}, we conclude that there exists $\sigma>0$ such that 
\begin{equation}\label{eq:vp great than u+ on T}
\tilde\vp>u^{+} +2\zeta \;\;\text{on}\;\;[T-\sigma,T[\;\times\;\mathbb{T}^{*}.
\end{equation}
More precisely, if \eqref{eq:vp great than u+ on T} does not hold for any $\sigma>0$, then there exists a sequence $(t_n,x_n)\in\D_{i}$ such that $t_n\rightarrow T$, $x_n\in \mathbb{T}^{*}$ and
$
\tilde\vp(t_n,x_n)\leq u^{+}(t_n,x_n) +2\zeta.
$
The compactness of $\mathbb{T}^{*}$ implies that there is a subsequence of $(t_n,x_n)$ which converges to $(T,x')$ for some $x'\in \mathbb{T}^{*}$. By taking the $\limsup$ of the above equation over the subsequence, we get
$
\tilde\vp(T,x')\leq u^{+}(T-,x')+2\zeta.
$
This contradicts \eqref{eq:vp great than u+ on DT}. Therefore, \eqref{eq:vp great than u+ on T} holds.  \\
\indent In \eqref{eq:vp great than u+ on T}, we choose $\sigma<\eps$. By a Dini-type argument, there exists $n_1\geq n_0$ such that  
\begin{equation}\label{eq:comp3_boundary_u+}
\tilde\vp>w_{n_1} +\zeta \;\;\text{on}\;\;[T-\sigma,T[\;\times\;\mathbb{T}^{*}.
\end{equation}
Set $w=w_{n_1}$. For $\kappa\in\; ]0,\varepsilon\wedge\alpha\wedge\zeta[\;$, define
$$w^{\kappa}:=
\left \{ 
\begin{array}{ll}
(\tilde\varphi -\kappa)\wedge w\ \ {\textrm on}\ \  [T-\sigma,T]\times \text{cl}(B_{\eps}(x_0)),\\
w \ \ \textrm{outside}\ \  [T-\sigma,T]\times \text{cl}(B_{\eps}(x_0)).
\end{array}
\right.
$$
Since $w(T,x)\geq g(x)$ and \eqref{eq: vp great than g} holds, we get that $w^{\kappa}(T,x)\geq g(x)$ for all $x\in\R^d$. We also notice that 
\begin{equation}\label{eq:boundary u+ contra}
w^{\kappa}(T,x_0)\leq \vp(x_0)-\kappa<u^{+}(T-,x_0)\leq u^{+}(T,x_{0}).
\end{equation} 
Using \eqref{eq: Lipschitz map exists_u+_bdd}, \eqref{eq: drift great than 0 boundary u+}, \eqref{eq:comp1_boundary_u+}, \eqref{eq:comp2_boundary_u+} and \eqref{eq:comp3_boundary_u+} in a manner that is similar to Step 1 in Theorem \ref{thm: main theorem_interior}'s proof, we can show that $w^{\kappa}$ is a stochastic super-solution, which contradicts \eqref{eq:boundary u+ contra}. \\
{\bf Step 2 (The super-solution property on $\DT$).} We will divide the proof into two steps: \\
{\bf Step 2.A.} We will show that $u^{-}(T-,\cdot)$ is a viscosity super-solution of 
$(\vp(x)-g(x))\mathbbm{1}_{\{H^{*}\vp(x)<\infty\}}\geq 0 \text{ on }\R^d.$
Let $x_0\in\R^d$ and $\vp\in C^2(\R^d)$ be such that $0=(u^{-}(T-,x_0)-\vp(x_0))=\min_{x\in\R^d}(u^{-}(T-,x)-\vp(x)).$
Assuming that $H^{*}\vp(x_0)=C<\infty$ and that $g(x_0)>u^{-}(T-,x_0)=\vp(x_0)$, we will work towards a contradiction. 
Let $\{w_{k}\}_{1}^{\infty}$ be a sequence in $\bU^{-}$ such that $ w_n\nearrow u^-$. Let
$ \tilde\vp(t,x)=\vp(x)-\iota|x-x_0|^{n_{0}}-(C+2)(T-t)$ and $\tilde\vp'(x)=\vp(x)-\iota|x-x_0|^{n_{0}}$ for $\iota>0$, where $\iota$ will be fixed later and $n_{0}\geq 2$ satisfies
 \begin{equation}\label{eq: test function decays faster_boundary u-}
\max_{0\leq t\leq T} (\tilde\vp(t,x)-w_{1}(t,x))\rightarrow - \infty \;\;\text{and}\;\; \max_{0\leq t\leq T} \tilde\vp(t,x)\rightarrow - \infty\;\; \text{as}\;\;|x|\rightarrow\infty \;\;\text{for any} \;\;\iota>0. 
 \end{equation}
 Note that $D\tilde\vp'(x)=D\tilde{\vp}(t,x)$ and $D^2\tilde\vp'(x)=D^2\tilde{\vp}(t,x)$. From $g(x_0)>\vp(x_0)=\tilde\vp(T,x_0)=u^{-}(T-,x_0)$ and the lower semi-continuity of $g$ and $u^-$, we can find $\eps>0$ and $\eta\in \; ]0,1[$ such that
 \begin{equation}\label{eq: vp less than u- in a neighborhood, u- boundary Step A}
g(x)-\tilde\vp(t,x)>\eps \;\;\text{for}\;\;(t,x)\in \text{cl}(B_{\eps}(T,x_0)),\;\;\tilde\vp<u^{-}+2\eta\;\;\text{on} \;\;[T-\eps,T[\;\times\; \text{cl}(B_{\eps/2}(x_0)).
\end{equation}
By the locally boundedness of $\mu_X$, $\sigma_X$, $\mu_Y$, $b$ and $\beta$, and $H^{*}\vp(x_0)=C$,  there exists $\iota>0$ such that
\begin{equation*}
\begin{array}{c}
\mu_Y(t,x,y,u)-\mu_{X}^{\top}(t,x,u) D\tilde\vp(t,x) -\frac{1}{2}\text{Tr}[\sigma_X\sigma_X^{\top}(t,x,u)D^2\tilde{\vp}(t,x)] \leq C+1\;\;\text{for all}\;\; (t,x,y, u)\in \D\times\R\times U
\end{array}
\end{equation*}
satisfying $(t,x)\in[T-\eps, T]\times \text{cl}(B_{\eps}(x_0)), |y-\tilde\vp(t,x)|\leq\eps$ and 
$u\in\mathcal{N}_{\eps, -\eta}(t,x,y,D\tilde{\vp}(t,x),\tilde\vp')$. Here, we may need to choose smaller values of $\eps$ and $\eta$. Therefore, by the definition of $\Delta^{u,e}$,
\begin{equation*}
\begin{array}{c}
\mu_Y(t,x,y,u)-\mathscr{L}^u\tilde{\vp}(t,x)\leq C+1-C-2 \leq -\eta \text{ for all }(t,x,y) \in \D\times\R\times U \\ \text{s.t. } (t,x)\in[T-\eps, T]\times \text{cl}(B_{\eps}(x_0)), |y-\tilde\vp(t,x)|\leq \eps \text{ and } u\in\mathcal{N}_{\eps,-\eta}(t,x,y,D\tilde{\vp}(t,x),\tilde{\vp}).
\end{array}
\end{equation*}
Fix $\iota$. Since $w_{k}\nearrow u^{-}$, there exists $n_0\in\N$ such that 
$
\tilde\vp<w_{n_0}+\eta\;\;\text{on} \;\;[T-\eps,T[\;\times\;\text{cl}(B_{\eps/2}(x_0))
$
due to \eqref{eq: vp less than u- in a neighborhood, u- boundary Step A}. By \eqref{eq: test function decays faster_boundary u-}, there exists $R_0>\eps$ such that
$\tilde\vp(t,x)<w_{n_0}(t,x)+\eps\leq w_n(t,x)+\eps\text{ on } \mathbb{O}$ for $n\geq n_{0}$, where $\mathbb{O} := [T-\eps,T]\times(\R^d\setminus\text{cl}(B_{R_0}(x_0)))$.
Since $\tilde{\vp}(T,x)\leq \vp(x)$, $u^{-}(T-,\cdot)-\tilde\vp(T,\cdot)$ is strictly positive on the compact set $\mathbb{T}^{*}:=\text{cl}(B_{R_0}(x_0))-B_{\eps/2}(x_0)$. Hence, by the lower semi-continuity of $u^{-}$, there exists $\alpha>0$ such that $
\tilde\vp(T,\cdot)<u^{-}(T-,\cdot)-4\alpha$ on $\mathbb{T}^{*}.$
Similar to Step 1 in this proof, we can find $\sigma\in \;]0,\eps[$ and $n_1\geq n_0$ such that
$
\tilde\vp<w_{n_1} -\alpha$ on $ [T-\sigma,T[\;\times\; \mathbb{T}^{*}.
$
Set $w=w_{n_1}$. For $\kappa \in\; ]0,\varepsilon\wedge\delta\wedge \alpha[\;$, define
$$w^{\kappa}:=
\left \{ 
\begin{array}{ll}
(\tilde\varphi +\kappa)\vee w\ \ {\textrm on}\ \  [T-\sigma,T]\times \text{cl}(B_{\eps}(x_0)),\\
w \ \ \textrm{outside}\ \  [T-\sigma,T]\times \text{cl}(B_{\eps}(x_0)).
\end{array}
\right.
$$
As in Step 2 of Theorem \ref{thm: main theorem_interior}'s proof, we can show that $w^{\kappa}\in \mathbb{bU}^{-}$, which yields a contradiction.  \\
\noindent {\bf Step 2.B:} In this step, we prove that $u^-(T-,\cdot)$ is a viscosity super-solution of $
\delta^{*}\vp(x) \geq 0$. Let $x_0\in\R^d$ and $\vp\in C^2(\R^d)$ be such that  
$
0=(u^{-}(T-,x_0)-\vp(x_0))=\min_{R^d}(u^{-}(T-,x)-\vp(x)).
$
Let $(s_n,\xi_n)$ be a sequence in $\D_{i}$ satisfying $(s_n,\xi_n)\rightarrow(T,x_0)$ and $u^{-}(s_n,\xi_n)\rightarrow u^{-}(T-,x_0)=\vp(x_0).$ For all $n\in \N$, $k\geq 0$ and $\iota\geq 0$, define
\begin{equation*}
\vp_{n}^{k,\iota}(t,x)=\vp(x)-\iota|x-x_0|^4+k\frac{T-t}{(T-s_n)}, \vp^{\iota}(x)=\vp(x)-\iota|x-x_0|^4.
\end{equation*}
Notice that 
\begin{equation*}\label{eq: limsup diff 0}
\lim_{\iota\rightarrow 0}\lim_{k \rightarrow 0}\limsup_{n\rightarrow \infty}\sup_{(t,x)\in[s_n, T]\times \text{cl}(B_1(x_0))}|\vp_n^{k,\iota}(t,x)-\vp(x)|=0.
\end{equation*}
Let $(t_n^{k,\iota}, x_n^{k,\iota})$ be the minimizer of $u^--\vp_n^{k,\iota}$ on $[s_n, T]\times \text{cl}(B_1(x_0))$. We claim that for any $k>0$ and $\iota>0$, there exists $N^{k,\iota}\in\N$ such that 
\begin{equation}\label{eq: convergence claim}
s_n\leq t_n^{k,\iota}<T \text{ for all } n\geq N^{k,\iota},\; \text{and}\;\; x_n^{k,\iota}\rightarrow x_0 \text{ as } n\rightarrow \infty.
\end{equation}
We now prove \eqref{eq: convergence claim}. Since $(s_n,\xi_n)\rightarrow(T,x_{0})$, we can find $N^{k,\iota}\in\N$ such that for $n\geq N^{k,\iota}$,
\begin{equation}\label{eq: t_n^k less than T 1}
(u^{-}-\vp_n^{k,\iota})(s_n,\xi_n)=u^{-}(s_n,\xi_n)-\vp(\xi_n)+\iota|\xi_n-x_0|^4-\frac{1}{k}\leq -\frac{1}{2k}<0.
\end{equation}
On the other hand, 
\begin{equation}\label{eq: t_n^k less than T 2}
\liminf_{t\uparrow T, x'\rightarrow x}(u^{-}-\vp_n^{k,\iota})(t,x')=u^{-}(T-,x)-\vp(x)+\iota|x-x_0|^4\geq 0\text{   for  }|x-x_0|\leq 1.
\end{equation}
By \eqref{eq: t_n^k less than T 1} and \eqref{eq: t_n^k less than T 2}, the first part of \eqref{eq: convergence claim} holds. By an argument similar to Step 4 in Theorem 3.1's proof in \cite{Bayraktar_and_Sirbu_SP_HJBEqn},  we know that the second part of \eqref{eq: convergence claim} also holds. \\
\indent From \eqref{eq: convergence claim} and the definition of $\vp_n^{k,\iota}$, we also see that
\begin{equation}\label{eq: convergence of vp_n at t_n x_n}
\vp_n^{k,\iota}(t_n^{k,\iota},x_n^{k,\iota})\rightarrow u^{-}(T-,x_0)=\vp(x_0)\;\;\text{as}\;\;n\rightarrow \infty, \text{ then  } k\rightarrow 0, \iota\rightarrow 0.
\end{equation}
By \eqref{eq: convergence claim}, \eqref{eq: convergence of vp_n at t_n x_n} and the facts that $u^-(t_n^{k,\iota},x_n^{k,\iota})\leq \vp_n^{k,\iota}(t_n^{k,\iota},x_n^{k,\iota}) \text{ and } \liminf_{(t<T, x)\rightarrow (T,x_0)}u^-(t,x)=u^-(T-,x_0),$
it holds that
$
u^-(t_n^{k,\iota},x_n^{k,\iota})\rightarrow u^{-}(T-,x_0)=\vp(x_0)$ as $n\rightarrow \infty$ then $k\rightarrow 0, \iota\rightarrow 0.
$
Since for all $k>0$, $\iota>0$ and $n\geq N^{k,\iota}$, $(t_n^{k,\iota}, x_n^{k,\iota})$ is a local minimizer of $u^--\vp_n^{k,\iota}$ and $t_n^{k,\iota}<T$, we get 
\begin{equation*}
-\partial_t\vp_n^{k,\iota}(t_n^{k,\iota},x_n^{k,\iota})+H^*(t_n^{k,\iota},x_n^{k,\iota},u^{-}(t_n^{k,\iota},x_n^{k,\iota}), D\vp_n^{k,\iota}(t_n^{k,\iota},x_n^{k,\iota}), D^2\vp_n^{k,\iota}(t_n^{k,\iota},x_n^{k,\iota}) )\geq 0
\end{equation*}
from Theorem \ref{thm: main theorem_interior}. For any $k>0$, $\iota>0$ and $n\geq N_n^{k,\iota}$, from the definition of $H^*$, there exists a sequence $\{(\eps_{m},\eta_{m},t_{m},x_{m},y_{m},p_{m},A_{m},\vp_{m})\} \subset \R_{+}\times [-1,1] \times\D\times\R\times \R^{d}\times\mathbb{M}^{d}\times C(\D)$ such that ($\eps_{m}, \eta_{m})\rightarrow(0,0)$,
\begin{gather}
\label{eq:convergence_u-_bbd}
\vp_m  \overset{\text{u.c.}}{\longrightarrow} \vp_n^{k,\iota}, 
(t_m, x_m, y_m, p_m, A_m)\rightarrow (t_n^{k,\iota},x_n^{k,\iota}, u^{-}(t_n^{k,\iota},x_n^{k,\iota}), D\vp_n^{k,\iota}(t_n^{k,\iota},x_n^{k,\iota}), D^2\vp_n^{k,\iota}(t_n^{k,\iota},x_n^{k,\iota}))\;\;\text{and} \\
H_{\eps_m,\eta_m}(t_m, x_m, y_m, p_m, A_m, \vp_m) \rightarrow H^*(t_n^{k,\iota},x_n^{k,\iota},u^{-}(t_n^{k,\iota},x_n^{k,\iota}), D\vp_n^{k,\iota}(t_n^{k,\iota},x_n^{k,\iota}), D^2\vp_n^{k,\iota}(t_n^{k,\iota},x_n^{k,\iota}) )>-\infty. \nonumber
\end{gather}
This implies that $\mathcal{N}_{\eps_m, \eta_m}(t_m, x_m, y_m, p_m,\vp_m)\neq \emptyset$ since $\sup \emptyset = -\infty$. By the definition of $\delta$, it holds that $\delta(t_m, x_m, y_m,p_m,\vp_m)\geq - \sqrt{\eps_m^2+\eta_m^2}.$
From \eqref{eq:convergence_u-_bbd} and the definition of $\delta^*$, we get
\begin{equation*}
\delta^{*}(t_n^{k,\iota},x_n^{k,\iota},u^{-}(t_n^{k,\iota},x_n^{k,\iota}), D\vp_n^{k,\iota}(t_n^{k,\iota},x_n^{k,\iota}), \vp_n^{k,\iota})\geq \limsup_{m\rightarrow\infty} \delta(t_m, x_m, y_m,p_m,\vp_m)\geq 0.
\end{equation*}
By the definition of $\Delta^{u,e}$ in the set-valued map $\mathbf{N}$, the equation above implies that
\begin{equation}\label{eq:new}
\delta^{*}(t_n^{k,\iota},x_n^{k,\iota}, u^-(t_n^{k,\iota}, x_n^{k,\iota}),  D\vp_n^{k,\iota}(t_n^{k,\iota},x_n^{k,\iota}),\vp^\iota)=\delta^{*}(t_n^{k,\iota},x_n^{k,\iota},u^{-}(t_n^{k,\iota},x_n^{k,\iota}), D\vp_n^{k,\iota}(t_n^{k,\iota},x_n^{k,\iota}), \vp_n^{k,\iota})\geq 0.
\end{equation}
Note that $\vp^\iota \overset{\text{u.c.}}{\longrightarrow} \vp$ as $\iota\rightarrow 0$. Moreover, for $\iota>0$, $u^-(t_n^{k,\iota},x_n^{k,\iota})\rightarrow\vp(x_0)$ and $D\vp_n^{k,\iota}(t_n^{k,\iota},x_n^{k,\iota})\rightarrow D\vp(x_0)$ as $n\rightarrow \infty$ then $k\rightarrow 0$. 
Taking the $\limsup$ of \eqref{eq:new} by first sending $n\rightarrow \infty$ then $k\rightarrow 0$ and $\iota\rightarrow 0$, we have $\delta^*\vp(x_0)=\delta^*\vp(T,x_0,\vp(x_0),D\vp(x_0),\vp)\geq 0$ from the upper semi-continuity of $\delta^*$, 
\qed
 
\section{Verification by Comparison}\label{sec: Verification}
\noindent We now carry out the verification for non-smooth functions assuming the comparison principle as in \cite{Bayraktar_and_Sirbu_SP_HJBEqn}. 
\begin{assumption}\label{assump: H continuous in domain}
Let $H=H_*$. Assume that $H=H^*$ on the set $\{H<\infty\}$ and that there exists an LSC function $G:\D\times\times\R\times\R^d\times\mathbb{M}^d\times C(\D)\rightarrow\R$ such that
\begin{equation*}
\begin{array}{c}
(a)\; H(t,x,y,p,A,\vp)<\infty \implies G(t,x,y,p,A,\vp)\leq 0, \\
(b)\;G(t,x,y,p,A,\vp)<0\implies H(t,x,y,p,A,\vp)<\infty.
\end{array}
\end{equation*} 
\end{assumption}
\begin{proposition}\label{prop: viscosity property after introducing G}
Under Assumptions \ref{assump: lambda intensity kernel}-\ref{assump: regularity} (resp. \ref{assump: lambda intensity kernel}-\ref{assump:semisolution_not_empty}) and \ref{assump: H continuous in domain}, $u^+$ (resp. $u^-$) is a USC (resp. an LSC) viscosity sub-solution (resp. super-solution) of $\max\left\{-\partial_t\vp(t,x)+H\vp(t,x),G\vp(t,x)\right\}=0\text{ on }\D_{i}.$
Moreover, if $g$ is USC, $u^{+}(T-,\cdot)$ is a USC viscosity sub-solution of 
$
\min\left\{\max\{\vp(x)-g(x),G\vp(x)\},\delta_{*}\vp(x)\right\}\leq 0\text{ on } \R^d.
$
If $g$ is LSC, $u^{-}(T-,\cdot)$ is an LSC viscosity super-solution of 
$
\min\left\{\max\{\vp(x)-g(x), G\vp(x)\},\delta^{*}\vp(x) \right\}\geq 0\text{ on }\R^d.
$
\end{proposition}
{\it Proof }
\textbf{(1) The sub-solution property in $\Di$.} Suppose
$0=(u^+-\varphi)(t_0,x_0)= \max_{ \Di}(u^+-\varphi)
$ for some $(t_0,x_0)\in\Di$ and $\varphi\in C^{1,2}(\D)$.
Then $-\partial_t\vp(t_0,x_0)+H\vp(t_0,x_0)=-\partial_t\vp(t_0,x_0)+H_*\vp(t_0,x_0)\leq 0$ from Theorem \ref{thm: main theorem_interior}.  
From (a) in Assumption \ref{assump: H continuous in domain}, $G\vp(t_0,x_0)\leq 0$. Therefore, the sub-solution property holds for $u^+$ in the parabolic interior. \\
\textbf{(2) The super-solution property in $\Di$.}  Suppose
$
0=(u^--\varphi)(t_0,x_0)= \min_{ \Di}(u^--\varphi)
$ for some $(t_0,x_0)\in\Di$ and $\varphi\in C^{1,2}(\D)$.
If $H\vp(t_0,x_0)<\infty$,
$-\partial_t\vp(t_0,x_0)+H\vp(t_0,x_0)=-\partial_t\vp(t_0,x_0)+H^*\vp(t_0,x_0)\geq 0$ from Assumption \ref{assump: H continuous in domain} and Theorem \ref{thm: main theorem_interior}. 
On the other hand, if $H\vp(t_0,x_0)=\infty$, $G\vp(t_0,x_0)\geq 0$ from (b) in Assumption \ref{assump: H continuous in domain}. Therefore, the viscosity super-solution property holds for $u^-$ in the parabolic interior. \\
\textbf{(3) The sub-solution property on $\DT$.}  From Theorem \ref{thm: bd_viscosity_property}, we know that $u^+(T-,\cdot)$ is viscosity sub-solution of $
\min\{\vp(x)-g(x), \delta_{*}\vp(x)\} \leq 0. $
Therefore, it suffices to show that $Gu^+(T-,\cdot)\leq 0$ in the viscosity sense. Let $x_0\in\R^d$ and $\vp\in C^2(\R^d)$ be such that  
$
0=(u^{+}(T-,x_0)-\vp(x_0))=\max_{x\in\R^d}(u^{+}(T-,x)-\vp(x)).
$
Let $(s_n,\xi_n)$ be a sequence in $\D_{i}$ satisfying $(s_n,\xi_n)\rightarrow(T,x_0) $ and $u^{+}(s_n,\xi_n)\rightarrow u^{+}(T-,x_0).$
For all $n\in \N$, $k\geq 0$ and $\iota\geq 0$, define
\begin{equation*}
\vp_{n}^{k,\iota}(t,x)=\vp(x)+\iota|x-x_0|^4-k\frac{T-t}{(T-s_n)}, \vp^{\iota}(x)=\vp(x)+\iota|x-x_0|^4.
\end{equation*}
Let $(t_n^{k,\iota}, x_n^{k,\iota})$ be the maximizer of $u^+-\vp_n^{k,\iota}$ on $[s_n, T]\times \text{cl}(B_1(x_0))$. Similar to the arguments in Step 2B of Theorem \ref{thm: bd_viscosity_property}'s proof, we can show that 
$
\lim_{k\rightarrow 0, \iota\rightarrow 0}\lim_{n\rightarrow \infty}  u^+(t_n^{k,\iota},x_n^{k,\iota}) = \vp(x_0).
$
We also know that for any $k>0$ and $\iota>0$, there exists $N^{k,\iota}\in\N$ such that $s_n\leq t_n^{k,\iota}<T \text{ for all } n\geq N^{k,\iota}$ and $x_n^{k,\iota}\rightarrow x_0 \text{ as } n\rightarrow \infty.$
Therefore, for all $k>0$, $\iota>0$ and $n\geq N^{k,\iota}$, $(t_n^{k,\iota}, x_n^{k,\iota})$ is a maximizer of $u^+-\vp_n^{k,\iota}$ on $[s_n, T]\times \text{cl}(B_1(x_0))$. From Theorem \ref{thm: main theorem_interior}, 
\begin{equation*}
-\partial_t\vp(t_n^{k,\iota},x_n^{k,\iota})+H_*(t_n^{k,\iota},x_n^{k,\iota}, u^+(t_n^{k,\iota},x_n^{k,\iota}), D\vp_n^{k,\iota}(t_n^{k,\iota},x_n^{k,\iota}), D^2\vp_n^{k,\iota}(t_n^{k,\iota},x_n^{k,\iota}),\vp_n^{k,\iota})\leq 0.
\end{equation*}
Hence, the $H_{*}$-term in the above equation is less than $\infty$. From the definition of $\Delta^{u,e}$, we get
\begin{equation*}
\begin{array}{c}
H_*(t_n^{k,\iota},x_n^{k,\iota}, u^+(t_n^{k,\iota},x_n^{k,\iota}), D\vp_n^{k,\iota}(t_n^{k,\iota},x_n^{k,\iota}), D^2\vp_n^{k,\iota}(t_n^{k,\iota},x_n^{k,\iota}), \vp^{\iota})<\infty, \;\;\text{which further implies that } \\ 
G\vp(t_n^{k,\iota},x_n^{k,\iota}, u^+(t_n^{k,\iota},x_n^{k,\iota}), D\vp_n^{k,\iota}(t_n^{k,\iota},x_n^{k,\iota}), D^2\vp_n^{k,\iota}(t_n^{k,\iota},x_n^{k,\iota}), \vp^{\iota})\leq 0\;\;\text{by Assumption \ref{assump: H continuous in domain}.}
\end{array}
\end{equation*}
 Using an argument similar to that in Step 2B of Theorem \ref{thm: bd_viscosity_property}'s proof, we conclude that $G\vp(x_0)\leq 0.$\\
\textbf{(4) The super-solution property on $\DT$.} It suffices to show that $u^-(T-,\cdot)$ is a viscosity super-solution of
\begin{equation}\label{eq:super_introduce_G}
\max\{\vp(x)-g(x), G\vp(x)\}\geq 0.
\end{equation} 
Let $x_0\in\R^d$ and $\vp\in C^2(\R^d)$ be such that  
$
0=(u^{-}(T-,x_0)-\vp(x_0))=\min_{x\in\R^d}(u^{-}(T-,x)-\vp(x)).
$
From Theorem \ref{thm: bd_viscosity_property}, one of the following two scenarios must hold:
\begin{gather}
\label{eq:scenario1_prop_introduce_G}
\vp(x_0)\geq g(x_0),\; H^*\vp(x_0)<\infty\;\;\text{or}\\
\label{eq:scenario2_prop_introduce_G}
H^*\vp(x_0)=\infty.
\end{gather}
\eqref{eq:scenario1_prop_introduce_G} implies \eqref{eq:super_introduce_G}; on the other hand, if \eqref{eq:scenario2_prop_introduce_G} holds, then
$H\vp(x_0)=\infty,$ which means that $G\vp(x_0)\geq 0$ from (b) in Assumption \ref{assump: H continuous in domain}. Therefore, \eqref{eq:super_introduce_G} holds. \qed
\begin{assumption}\label{ass: cmp_1}
Assume that $\delta^*=\delta_*$, $g$ is continuous and a comparison principle holds between USC sub-solutions and LSC super-solutions for
\begin{equation}\label{eq: bd_pde}
\min\{\max\{\vp(x)-g(x), G\vp(x)\}, \delta\vp(x)\}=0\;\;\text{on}\;\;\R^{d}.
\end{equation}
\end{assumption}
In the presence of jumps, it is nontrivial to check this assumption. When there are no jumps in the controlled processes, the comparison principle can be proved in certain
classes of functions (see the discussion above Assumption 2.2 in \cite{Bouchard_Elie_Touzi_ControlledLoss}). Also, in Section \ref{sec:equivalence}, $\delta$ drops out in the corresponding PDE and there are comparison results available for fully non-linear equations with jumps (see \cite{Barle-IntegroPDE}). 
\begin{lemma}\label{lem: bd_equal}
Under Assumptions \ref{assump: H continuous in domain}, \ref{ass: cmp_1}
and \ref{assump: lambda intensity kernel}-\ref{assump: regularity},
$
u^{-}(T-,\cdot)=u^{+}(T-,\cdot)=\hat{g}(\cdot),
$
where $\hat{g}$ is the unique continuous viscosity solution to \eqref{eq: bd_pde}.
\end{lemma}
{\it Proof } It follows from their definitions that $u^- \leq u^+$. Since $u^+$ is USC and $u^{-}$ is LSC, then
\begin{equation*}
\label{eq:limits}
u^-(T-,x)=\liminf _{(t<T, x')\rightarrow (T,x)} u^-(t,x')\leq \limsup _{(t<T, x')\rightarrow (T,x)}u^+(t,x')= u^+(T-,x).
\end{equation*}
Moreover, $u^{+}(T-,\cdot)$ is a viscosity sub-solution and $u^{-}(T-,\cdot)$ is a viscosity super-solution to \eqref{eq: bd_pde} due to Theorem \ref{thm: bd_viscosity_property}. Therefore, the claim holds by Assumption \ref{ass: cmp_1}.
\qed
\begin{theorem}\label{thm: unique viscosity solution}
Suppose that there is a comparison principle for 
\begin{equation}\label{eq:HJB}
\max\{-\partial_t\vp(t,x)+H\vp(t,x),\; G\vp(t,x)\}=0 \text{ on }\Di
\end{equation}
and that Assumptions \ref{assump: lambda intensity kernel}-\ref{assump: regularity}, \ref{assump: H continuous in domain} and \ref{ass: cmp_1} hold. Then there exists a unique continuous viscosity solution
$V$ to \eqref{eq:HJB} with terminal condition $V(T,\cdot)=\hat{g}(\cdot)$ and
$u(t,x)=u^{-}(t,x)=u^{+}(t,x)=V(t,x)$ for $(t,x)\in \Di.$
\end{theorem}
{\it Proof } Define
$$\hat{u}^{+}(t,x):=
\left \{ 
\begin{array}{ll}
u^{+}(t,x), \;\;& (t,x)\in\Di\\
\hat{g}(x),\;\; & t=T, x\in\R^d
\end{array}
\right.
\;\;\text{and}\;\;
\hat{u}^{-}(t,x):=
\left \{ 
\begin{array}{ll}
u^{-}(t,x), \;\;& (t,x)\in\Di,\\
\hat{g}(x),\;\; & t=T, x\in\R^d.
\end{array}
\right.
$$
From Proposition \ref{prop: viscosity property after introducing G}, $\hat{u}^-$ is an LSC viscosity super-solution
and $\hat{u}^+$ is a USC viscosity sub-solution of \eqref{eq:HJB}. Since $
\hat{u}^{+}(T,\cdot) = \hat{u}^{-}(T, \cdot)$,  $\hat{u}^{+}\leq \hat{u}^-$ on $\D$ by comparison. Hence, $\hat{u}^{+} = \hat{u}^{-}$ on $\D$ from \eqref{eq:intfvmavp}. Define $V:=\hat{u}^{+} = \hat{u}^{-}$. It is a continuous viscosity solution of \eqref{eq:HJB} satisfying $V(T,x)=\hat{g}(x)$. Uniqueness follows directly from the comparison principle. 
\qed

\section{Stochastic Control as a Stochastic Target Problem}\label{sec:equivalence}
\noindent In this section, we show how the HJB equation associated to an optimal control problem in standard form can be deduced from a stochastic target problem. Given a bounded continuous function $g:\R^d\rightarrow\R$, we define an optimal control problem by
$
 \mathbf{u}(t,x):=\inf_{\nu\in\mathcal{U}^t}\mathbb{E}[g(X_{t,x}^{\nu}(T))]. 
$
We follow the setup of Section \ref{sec:prob} with one exception:
$\mathcal{U}^{t}$ is the collection of all the $\mathbb{F}^{t}$-predictable processes in $\mathbb{L}^2(\Omega\times[0,T], \mathcal{F}\otimes\mathcal{B}[0,T], \mathbb{P}\otimes \lambda_{L}; U)$, where $U\subset \R^{d}$ and $X$ follows the SDE
\begin{equation*}
dX(s)=\mu_{X}(s,X(s),\nu(s))ds+\sigma_{X}(s,X(s),\nu(s))dW_s+\int_{E} \beta(s,X(s-),\nu(s), e)\lambda(ds,de). 
\end{equation*}
To convert the control problem to its stochastic target counterpart,  we need the following lemma,  which is an adaptation of a result in \cite{Bouchard_Equivalence}. 

\begin{lemma}\label{eq:equivalence_application}
Suppose Assumptions \ref{assump: lambda intensity kernel} and \ref{assump: regu_on_coeff} hold. Define a stochastic target problem as follows:
\begin{equation*}
\begin{gathered}
u(t,x):=\inf\{y\in\R: \exists (\nu, \alpha, \gamma)\in \mathcal{U}^t\times\mathcal{A}^t\times\Gamma^t\;\text{s.t.}\;Y_{t,y}^{\alpha,\gamma}(T)\geq g(X_{t,x}^{\nu}(T))\}, \text{where} \\
Y_{t,y}^{\alpha,\gamma}(\cdot):=y+\int_t^{\cdot}\alpha^{\top}(s)dW_s+\int_t^{\cdot}\int_E\gamma^{\top}(s,e)\tilde{\lambda}(ds,de) 
\end{gathered}
\end{equation*}
and $\mathcal{A}^t$ and $\Gamma^{t}$ are the collections of $\R^{d}$-valued and $\mathbb{L}^{2}(E, \mathcal{E}, \hat{m}; \R^{I})$-valued processes, respectively,  satisfying the admissibility conditions in Section \ref{sec:prob}. 
Then $u=\mathbf{u}$ on $\D$.
\end{lemma}
{\it Proof }
Since $\mathcal{A}^{t}$ and $\Gamma^{t}$ satisfy the admissibility conditions,  this stochastic target problem is well defined.  In view of Lemma 2.1 in \cite{Bouchard_Equivalence},  it suffices to check that
\begin{equation}\label{eq:inclusion_equi}
\left\{ g(X_{t,x}^{\nu}(T),\nu\in\mathcal{U}^{t}\right\}\subset \left\{ M(T), M\in \mathcal{M}\right\},\text{ where } 
\mathcal{M}:=\left\{Y_{t,y}^{\alpha,\gamma}(\cdot): y\in\R, \alpha\in\mathcal{A}^{t}, \gamma\in\Gamma^{t} \right\}.
\end{equation}
In fact, by the martingale representation theorem, for  any $\nu\in\mathcal{U}^{t}$, $\mathbb{E}[g(X_{t,x}^{\nu}(T))|\mathcal{F}^{t}_{\cdot}]$ can be represented in the form of $Y_{t,y}^{\alpha, \gamma}$ for some $\alpha \in\mathcal{A}^{t}$ and $\gamma\in\Gamma^{t}_{0}$, where  $\Gamma_{0}^{t}$ is the collection of $\mathbb{L}^{2}(E, \mathcal{E}, \hat{m}; \R^{I})$-valued processes satisfying all of the admissibility conditions except for $\eqref{eq:admissibility}$. In particular,  $g(X_{t,x}^{\nu}(T))=Y_{t,y}^{\alpha, \gamma}(T)$. 
Assume, contrary to \eqref{eq:inclusion_equi}, that there exists $\nu_{0}\in\mathcal{U}^{t}$ such that
\begin{equation*}
\mathbb{E}[g(X_{t,x}^{\nu_{0}}(T))|\mathcal{F}^{t}_{\cdot}]=y+\int_t^{\cdot}\alpha_{0}^{\top}(s)dW_s+\int_t^{\cdot}\int_E\gamma_{0}^{\top}(s,e)\tilde{\lambda}(ds,de)
\end{equation*}
for some $y\in\mathbb{R}$,  $\alpha_{0}\in\mathcal{A}^{t}$ and $\gamma_{0}\in\Gamma^{t}_{0}$, but  \eqref{eq:admissibility} does not hold. In the equation above, $\mathbb{E}[g(X_{t,x}^{\nu_{0}}(T))|\mathcal{F}^{t}_{\cdot}]$ can be chosen to be c\`adl\`ag, thanks to Theorem 1.3.13 in \cite{ShreveKaratzas}. Then for $K>2\|g\|_{\infty}$, there exists $\tau_{0}\in\mathcal{T}_{t}$ such that 
$\mathbb{P}\left(\left|\int_{E}\gamma^{\top}(\tau_{0},e) \lambda(\{\tau_{0}\},de)\right|>K\right)>0. $
Suppose that
$\mathbb{P}\left(\int_{E}\gamma^{\top}(\tau_{0},e) \lambda(\{\tau_{0}\},de)>K\right)>0. \footnote{If this does not hold, the integral is less than $-K$ with positive probability. Noticing this, we can carry out the proof in a similar manner when this assumption does not hold.}$
Let $M_{0}(\cdot)=\mathbb{E}\left[g(X_{t,x}^{\nu_{0}}(T))|\mathcal{F}^{t}_{\cdot}\right]$. Therefore, $$M_{0}(\tau_{0})-M_{0}(\tau_{0}-) = \int_{E}\gamma^{\top}(\tau_{0},e) \lambda(\{\tau_{0}\},de)>K\;\;\text{with positive probability}.$$
Since $|M_{0}|$ is bounded by $\|g\|_{\infty}<K/2$, we obtain a contradiction. \qed
Let $\mathbf{H}^*$ be the USC envelope of the LSC map $\mathbf{H}:\D\times\R^d\times\mathbb{M}^d\times C(\D) \rightarrow \R$ defined by
\begin{equation*}
\begin{array}{c}
\mathbf{H}: (t,x,p,A, \vp) \rightarrow \sup_{u\in U}\{-I[\vp](t,x,u)-\mu_X^{\top}(t,x,u)p-\frac{1}{2}\text{Tr}[\sigma_X\sigma_X^{\top}(t,x,u)A]\},\;\text{where} \\
I[\vp](t,x,u)=\sum_{1\leq i\leq I}\int_E \left( \vp(t,x+\beta_i(t,x,u,e))-\vp(t,x)\right)m_i(de).
\end{array}
\end{equation*}
\begin{theorem}\label{thm: optimal control}
Under Assumptions \ref{assump: lambda intensity kernel} and \ref{assump: regu_on_coeff}, $u^+$  is a USC viscosity sub-solution of 
$$
-\partial_t\vp(t,x)+\mathbf{H}\vp(t,x)\leq 0\text{ on } \Di
$$
and $u^+(T-,x)\leq g(x)$ for all $x\in \R^d$. On the other hand, $u^-$  is an LSC viscosity super-solution of 
$$
-\partial_t\vp(t,x)+\mathbf{H}^{*}\vp(t,x)\geq 0\text{ on } \Di
$$
and $u^-(T-,\cdot)$ is an LSC viscosity super-solution of 
$$
\left(\vp(x)-g(x)\right)\mathbbm{1}_{\{\mathbf{H}^*\vp(x) <\infty\}} \geq 0\text{ on }\R^d. 
$$
\end{theorem}
{\it Proof }
It is easy to check Assumption \ref{assump: regularity} for the stochastic target problem. Since $g$ is bounded, we can check that all of the assumptions in the Appendix A are satisfied, which implies that Assumption \ref{assump:semisolution_not_empty} holds. From Theorem \ref{thm: main theorem_interior}, 
$u^+$ is a USC viscosity sub-solution of 
$
-\partial_t\vp(t,x)+H_{*}\vp(t,x)\leq 0\text{ on }\Di
$
and $u^-$ is an LSC viscosity super-solution of 
$
-\partial_t\vp(t,x)+H^{*}\vp(t,x)\geq 0\text{ on }\Di.
$
From Proposition 3.1 in \cite{Bouchard_Equivalence}, 
$H^*\leq \mathbf{H}^*$ and $H_*\geq \mathbf{H}$. This implies that the viscosity properties in the parabolic interior hold. 

Also, by Theorem \ref{thm: bd_viscosity_property},
$u^{+}(T-,\cdot)$ is a USC viscosity sub-solution of 
$
\min\{\vp(x)-g(x), \delta_{*}\vp(x)\} \leq 0\text{ on }\R^d
$
and $u^-(T-,\cdot)$ is an LSC viscosity super-solution of 
$
\min\{\left(\vp(x)-g(x)\right)\mathbbm{1}_{\{H^*\vp(x) <\infty\}}, \delta^{*}\vp(x)\} \geq 0\text{ on }\R^d, 
$
where $\delta=\text{dist}(0,\mathbf{N}^c)-\text{dist}(0, \mathbf{N})$ and
\begin{equation*}
\begin{array}{ll}
\mathbf{N}(t,x,y,p,\vp)=&\{(q,s)\in\R^d\times\R: \exists (u,a,r)\in U\times\R^d\times\L^2(E,\mathcal{E}, \hat{m};\R^I)\;\text{s.t.}
\;q=a-\sigma_X^{\top}(t,x,u)p \\ &\text{and }s\leq \min_{1\leq i\leq I}\{r_i(e)-\vp(t,x+\beta_i(t,x,u,e))+\vp(t,x)\} \; \hat{m}-\text{a.s.}\; e\in E\; \}.
\end{array}
\end{equation*}
Obviously, $\mathbf{N}=\R^d\times\R$. Therefore, $\delta=\infty$ and the boundary conditions hold.
\qed
\noindent The following two corollaries show that $\mathbf{u}$ is the unique viscosity solution to its associated HJB equation. We omit the proof, since it is the same as the proofs of Proposition \ref{prop: viscosity property after introducing G} and Theorem \ref{thm: unique viscosity solution}.
\begin{corollary}\label{coro2}
Suppose that Assumptions \ref{assump: lambda intensity kernel} and \ref{assump: regu_on_coeff} hold,  $\mathbf{H}=\mathbf{H}^*$ on $\{\mathbf{H}<\infty\}$ and there exists an LSC function $\mathbf{G}:\D\times\R\times\R^d\times\mathbb{M}^d\times C(\D)\rightarrow\R$ such that 
\begin{equation*}
\begin{array}{c}
(a)\; \mathbf{H}(t,x,y,p,M,\vp)<\infty \implies \mathbf{G}(t,x,y,p,M,\vp)\leq 0, \\
(b)\;\mathbf{G}(t,x,y,p,M,\vp)<0\implies \mathbf{H}(t,x,y,p,M,\vp)<\infty.
\end{array}
\end{equation*} 
Then $u^+$ $(\text{resp. } u^-)$ is a USC $($resp. an LSC$)$ viscosity sub-solution $($resp. super-solution$)$ of
\begin{equation*}
\max\{-\partial_t\vp(t,x)+\mathbf{H}\vp(t,x), \mathbf{G}\vp(t,x)\}=0\;\;\text{on}\;\;\D_{i}
\end{equation*}  
and $u^{+}(T-,\cdot)$ $(\text{resp. } u^{-}(T-,\cdot))$ is a USC $($resp. an LSC$)$ viscosity sub-solution $($resp. super-solution$)$ of 
\begin{equation*}
\max\{\vp(x)-g(x), \mathbf{G}\vp(x)\}= 0\text{ on } \R^d.
\end{equation*}
\end{corollary}

\begin{corollary}
Suppose that all of the assumptions in Corollary \ref{coro2} hold.  Additionally, assume that there is a comparison principle between USC sub-solutions and LSC super-solutions for the PDE
\begin{equation}\label{eq:bd_pde_control}
\max\{\vp(x)-g(x), \mathbf{G}\vp(x)\}=0\;\;\text{on}\;\;\R^{d}.
\end{equation}
Then 
$u^{+}(T-,x)=u^{-}(T-,x)=\hat{\mathbf{g}}(x),$
where $\hat{\mathbf{g}}$ is the unique viscosity solution to \eqref{eq:bd_pde_control}. Furthermore, if the comparison principle holds for 
\begin{equation}\label{eq:HJB control}
\max\{-\partial_t\vp(t,x)+\mathbf{H}\vp(t,x),\; \mathbf{G}\vp(t,x)\}=0 \text{ on }\Di,
\end{equation}
then there exists a unique continuous viscosity solution
$\mathbf{V}$ to \eqref{eq:HJB control} with terminal condition $\mathbf{V}(T,x)=\hat{\mathbf{g}}(x)$ and
$
\mathbf{u}(t,x)=u(t,x)=u^{+}(t,x)=u^{-}(t,x)=\mathbf{V}(t,x)\;\;\text{for}\;\;(t,x)\in \Di.
$
\end{corollary}

\section{Conclusions}
In this paper, stochastic target problems in a jump diffusion setup are analyzed by using stochastic Perron's method, which had been recently developed to analyze the classical stochastic control problems. In fact, we using the fact that ordinary stochastic control problems can be embedded into stochastic target problems we extended that analysis to cover to processes in which both the diffusions and jumps are controlled.  
Our future research will focus on extending the analysis to stochastic target games. In the formulation of such problems, a strategic player tries to find a strategy such that the controlled process reaches a given target no matter what the opponent's control is. Of particular importance is the set-up in which one of the players is a stopper, whose aim is to get to the target at a stopping time instead of a fixed horizon. 

\begin{singlespace}
\appendix
\section*{Appendix A}\label{sec:appendix}
\renewcommand{\thesection}{A}
\renewcommand{\thetheorem}{A.\arabic{theorem}}
\noindent  We provide sufficient conditions for the nonemptiness of $\bU^+$ and $\bU^-$. 
\begin{assumption}\label{assump: g bounded}
$g$ is bounded.
\end{assumption}
 \begin{assumption} \label{assump: existence of no_investing_strategy}
 There exists $u_0 \in U$ such that $\sigma_Y(t,x,y,u_0)=0$ and $b(t,x,y,u_0(e),e)=0$ for all $(t,x,y,e)\in \D\times \mathbb{R}\times E$.
 \end{assumption}
 \begin{re}
 \textnormal{In the context of super-hedging in mathematical finance, the assumption above is equivalent to restricting trading to the riskless assets.}
 \end{re}
 \begin{proposition}\label{prop: U^+ is not empty}
 Under Assumptions  \ref{assump: lambda intensity kernel}, \ref{assump: regu_on_coeff}, \ref{assump: g bounded} and \ref{assump: existence of no_investing_strategy}, $\bU^{+}$ is not empty.
 \end{proposition}
{\it Proof }
\noindent\textbf{Step 1.}
In this step we assume that $\mu_{Y}$ is non-decreasing in its $y$-variable. We will show that $w(t,x)=\gamma-e^{kt}$ is a stochastic super-solution for some choice of $k$ and $\gamma$. 

By the linear growth condition on $\mu_Y$ in Assumption \ref{assump: regu_on_coeff}, there exists $L>0$ such that $|\mu_Y(t,x,y,u_0)|\leq  L(1+|y|)$, where $u_0$ is the element in $U$ in Assumption \ref{assump: existence of no_investing_strategy}. Choose $k\geq 2L$ and $\gamma$ such that $-e^{k T}+\gamma\geq \|g\|_{\infty}$. 
 Then $w(T,x)\geq g(x)$. It suffices to show that for any $(t,x,y)\in \D\times\mathbb{R}$, $\tau\in\mathcal{T}_t$, $\nu\in \Uc^t$ and $\rho\in \mathcal{T}_{\tau}$,
 \begin{equation}\label{eq: property of stochastic super-solution}
  Y(\rho )\geq w(\rho, X(\rho )) \;\; \mathbb{P}\text{-a.s.}\;\; \text{on}\;\;\{Y(\tau )\geq w(\tau, X(\tau))\}, \text{where } X:= X_{t,x}^{\nu\otimes_{\tau}u_0}, Y:=Y_{t,x,y}^{\nu\otimes_{\tau}u_0}.
 \end{equation}
Let $A= \{Y(\tau )> w(\tau, X(\tau))\}$, $V(s)=w(s,X(s))$ and $\Gamma(s)=\left(V(s)-Y(s)\right)\mathbbm{1}_{A}.$ Therefore, for $s\geq \tau$, 
\begin{gather}\label{eq: Gamma_integral}
dY(s)= \mu_{Y}\left(s,X(s),Y(s), u_0\right)ds, \;dV(s)= -ke^{ks}ds, \;
\Gamma(s)=\mathbbm{1}_{A}\int_{\tau}^{s} ( \xi(q)+ \Delta(q) ) dq, \text{where} \\
\Delta(s):=-ke^{ks}-\mu_Y(s,X(s),Y(s),u_0)\leq -ke^{ks}-\mu_Y(s,X(s),-e^{ks},u_0)\leq -ke^{ks}+L(1+e^{ks})\leq 0, \nonumber\\
\xi(s):=\mu_{Y}(s,X(s),V(s),u_0) -\mu_{Y}(s,X(s),Y(s),u_0). \nonumber
\end{gather}
Therefore, from \eqref{eq: Gamma_integral} it holds that
$$
\Gamma(s)\leq \mathbbm{1}_A\int_{\tau}^{s} \xi(q) dq \;\; \text{and}\;\; \Gamma^{+}(s)\leq \mathbbm{1}_A\int_{\tau}^{s} \xi^{+}(q) dq 
 \;\; \text{for} \;\; s\geq \tau.$$ 
From the Lipschitz continuity of $\mu_Y$ in $y$-variable in Assumption \ref{assump: regu_on_coeff}, 
\begin{equation*}\label{eq: fun_gronwall}
\Gamma^{+}(s)\leq \mathbbm{1}_A \int_{\tau}^{s} \xi^{+}(q) dq \leq \int_{\tau}^{s} L_0 \Gamma^{+}(q) dq  \;\; \text{for} \;\; s\geq \tau,
\end{equation*} 
where $L_0$ is the Lipschitz constant of $\mu_Y$ with respect to $y$. Note that we use the assumption that $\mu_Y$ is non-decreasing in its $y$-variable to obtain the second inequality.
Since $\Gamma^+(\tau)=0$, an application of Gr\"{o}nwall's Inequality implies that
 $\Gamma^+(\rho)\leq 0$, which further implies that \eqref{eq: property of stochastic super-solution} holds. \\
\textbf{Step 2.} We get rid of our assumption on $\mu_{Y}$ from Step 1 by following a proof similar to those in \cite{BayraktarLi} and \cite{Bouchard_Nutz_TargetGames}. For $c>0$, define $\widetilde Y_{t,x,y}^{\nu}$ as the strong solution of
\begin{equation*}
  d\widetilde{Y}(s)=\tilde \mu_{Y}(s,X_{t,x}^{\nu}(s),\widetilde{Y}(s),\nu(s)) ds +\tilde \sigma_{Y}^{\top}(s,X_{t,x}^{\nu}(s),\widetilde{Y}(s),\nu(s))dW_{s} + \int_{E} \widetilde{b}^{\top}(s,X_{t,x}^{\nu}(s-),\widetilde{Y}(s-), \nu_1(s),\nu_2(s,e), e)\lambda(ds,de)
\end{equation*}
  with initial data $\widetilde{Y}(t)=y$, where
  \begin{equation*}
   \widetilde{\mu}_{Y}(t,x,y,u):= c y+e^{ct} \mu_{Y}(t,x,e^{-c t}y,u), \;
  \widetilde{\sigma}_{Y}(t,x,y,u):= e^{c t} \sigma_{Y}(t,x,e^{-c t} y,u), \;
  \widetilde{b}(t,x,y,u(e),e):= e^{c t} b(t,x,e^{-c t} y,u(e),e). 
  \end{equation*}
 Therefore,  $$\widetilde{Y}_{t,x,y}^{\nu}(s)e^{-cs}=Y_{t,x,ye^{-ct}}^{\nu}(s), \;t\leq s\leq T. $$  Let 
  $
  \tilde u(t,x)= \inf\{y\in \R: \exists \; \nu\in \mathcal{U}^t, \mbox{  s.t.}\; \widetilde{Y}^{\nu}_{t,x,y}(T)\ge \tilde g(X^{\nu}_{t,x}(T))\;\as\},
  $
  where $\tilde{g}(x)=e^{c T} g(x)$. Therefore, $\tilde{u}(t,x)=e^{ct}u(t,x).$ Since $\mu_{Y}$ is Lipschitz in $y$, we can choose $c>0$ so that
  \begin{equation*}\label{eq: muY_uhat_tilde}
 \widetilde {\mu}_{Y}: (t,x,y,u) \mapsto  cy +  e^{c t}\mu_{Y}(t,x,e^{-c t}y,u)
  \end{equation*}
  is non-decreasing in $y$. Moreover, all the properties of $\widetilde{\mu}_{Y}, \widetilde{\sigma}_{Y}$ and $\widetilde{b}$ in Assumption \ref{assump: regu_on_coeff} still hold. We replace $\mu_Y$, $\sigma_Y$ and $b$ in all of the equations and definitions in Section \ref{sec:prob} with  $\widetilde{\mu}_{Y}, \widetilde{\sigma}_{Y}$ and $\widetilde{b}$, we get $\widetilde{H}^*$ and $\widetilde{H}_*$. Let $\widetilde{\bU}^+$ be the set of stochastic super-solutions of $$ -\partial_t\vp(t,x)+\widetilde{H}^*\vp(t,x)\geq 0\;\;\text{on}\;\;\Di.$$ It is easy to see that $w\in\bU^+$ if and only if $\widetilde{w}(t,x):=e^{ct}w(t,x)\in \widetilde{\bU}^+$. From Step 1, $\widetilde{\bU}^+$ is not empty. Thus, $\bU^+$ is not empty. 
\qed
 \begin{assumption}\label{assump: linear growth in y}
 There is $C\in \R$ such that for all $(t,x,y,u,e)\in\D\times\R\times U\times E$,
 $$\left|\mu_Y(t,x,y,u)+\int_E b^{\top}(t,x,y,u(e),e) m(de)\right|\leq C(1+|y|).$$
 \end{assumption}
\begin{proposition}\label{Prop: U^- is not empty}
Under Assumptions \ref{assump: lambda intensity kernel}, \ref{assump: regu_on_coeff}, \ref{assump: g bounded} and \ref{assump: linear growth in y}, $\mathbb{U}^{-}$ is not empty.
\end{proposition}
 {\it Proof }Assume that $$\mu_{Y}(t,x,y,u)+\int_E b^{\top}(t,x,y,u(e),e)m(de)$$ is non-decreasing in its $y$-variable. We could remove this assumption by using the argument from previous proposition. 

Choose $k\geq 2C$ ($C$ is the constant in Assumption \ref{assump: linear growth in y}) and $\gamma>0$ such that $e^{k T}-\gamma<- \|g\|_{\infty}$. Let $w(t,x)=e^{kx}-\gamma$. Notice that $w$ is continuous, has polynomial growth in $x$ and $w(T,\cdot)\leq g(\cdot)$. It suffices to show that for any $(t,x,y)\in \D\times\mathbb{R}$, $\tau\in\mathcal{T}_{t}$ and $\nu\in \Uc^t$, it holds that
$\mathbb{P}(Y(\rho )< w(\rho, X(\rho ))|B)>0$
for all $\rho \in \mathcal{T}_{\tau}$ and $B\subset \{Y(\tau)<w(\tau,X(\tau))\}$ satisfying $B\in\mathcal{F}_\tau^t$ and $\P(B)>0$, where $X:= X_{t,x}^{\nu}$ and $Y:=Y_{t,x,y}^{\nu}$. Define
 \begin{equation*}
 \begin{gathered}
M(\cdot)=Y(\cdot)-\int_{\tau}^{\cdot}K(s)ds,\;\; V(s)=w(s,X(s)),\;\;
 A= \{Y(\tau )<w(\tau, X(\tau))\},\;\; \Gamma(s)=\left(Y(s)-V(s)\right)\mathbbm{1}_{A}, \;\text{where} \\
K(s):=\mu_{Y}(s,X(s),Y(s),\nu(s))+\int_{E}b^{\top}(s,X(s-),Y(s-),\nu_1(s),\nu_2(s,e),e)m(de),\\
 \widetilde{K}(s):=\mu_{Y}(s,X(s),V(s),\nu(s))+\int_{E}b^{\top}(s,X(s-),V(s-),\nu_1(s),\nu_2(s,e),e)m(de).
 \end{gathered}
 \end{equation*}
 It is easy to see that $M$ is a martingale after $\tau.$  Due to the facts that $A\in\mathcal{F}_\tau^t$ and $dV(s)= ke^{ks}ds$, we further know
 \begin{equation}\label{eq: supermar1_nonemptyness of U+}
  \mathbbm{1}_{A}\left(Y(\cdot)-V(\cdot)+\int_{\tau}^{\cdot} ke^{ks}-K(s) ds \right)\;\; \text{is a super-martingale after}\;\;\tau .
 \end{equation} 
Since Assumption \ref{assump: linear growth in y} holds and $\mu_{Y}(t,x,y,u)+\int_E b^{\top}(t,x,y,u(e),e)m(de)$ is non-decreasing in $y$,
$$
\widetilde{K}(s)\leq \mu_Y(s,X(s),e^{ks}, \nu(s))+\int_{E}b^{\top}(s,X(s-),e^{ks},\nu_1(s),\nu_2(s,e),e)m(de)\leq 2C e^{ks}. 
$$
Therefore, it follows from \eqref{eq: supermar1_nonemptyness of U+} and the inequality above that 
\begin{equation}\label{eq: supermar2_nonemptyness of U+}
\widetilde{M}(\cdot):=\mathbbm{1}_{A}\left(Y(\cdot)-V(\cdot)-\int_{\tau}^{\cdot}\xi(s)ds)\right) \;\;\text{is a super-martingale after }\tau, \text{ where } \xi(s):=K(s)-\widetilde{K}(s).
\end{equation}
Since $\widetilde{M}(\tau)<0$ on $B$, there exists a non-null set $F\subset B $ such that $\widetilde{M}(\rho)<0$ on $F$. By the definition of $\widetilde{M}$ in \eqref{eq: supermar2_nonemptyness of U+}, we get
\begin{equation}\label{eq: Gamma_rho_strict_ineq on F}
\Gamma(\rho)< \mathbbm{1}_{A}\int_{\tau}^{\rho}\xi(s)ds \;\;\text{on}\;\;F.
\end{equation}
 Therefore, 
 \begin{equation}\label{eq: ineq for gronwall ineq}
 \Gamma^{+}(\rho)\leq \mathbbm{1}_A\int_{\tau}^{\rho} \xi^{+}(s) ds 
\leq \int_{\tau}^{\rho} L_0 \Gamma^{+}(s) ds\;\;\text{on}\;\;F. 
\end{equation}
  By Gr\"{o}nwall's Inequality,  $\Gamma^+(\tau)=0$ implies that
  $\Gamma^+(\rho)=0$ on $F$. More precisely, for $\omega\in F$ ($\P-\text{a.s.}$), $\Gamma^{+}(s)(\omega)=0$ for $s\in [\tau(\omega),\rho(\omega)]$. This implies that we can replace the inequalities with equalities in \eqref{eq: ineq for gronwall ineq}. Therefore, by \eqref{eq: Gamma_rho_strict_ineq on F}, $\Gamma(\rho)<0$ on $F$, which yields  
 $\mathbb{P}(Y(\rho )< w(\rho, X(\rho ))|B)>0.$ 
\qed

\section*{Appendix B}\label{sec:appendixB}
{\it Proof of Theorem \ref{thm: main theorem_interior}} \\
{\bf Step 1 ($u^+$ is a viscosity sub-solution).}
Assume, on the contrary, that for some $(t_0,x_0)\in\Di$ and $\varphi\in C^{1,2}(\D)$ satisfying 
$0=(u^+-\varphi)(t_0,x_0)= \max_{ \Di}(u^+-\varphi)$,
we have
\begin{equation}\label{eq: initial_local_contra_condition_at_t0x0_for u+}
4\eta:=-\partial_t\vp(t_0,x_0)+H_*\vp(t_0,x_0)>0.
\end{equation}
From Lemma \ref{lem: monotone seq approaches v+ or v_-}, there exists a non-increasing sequence $\mathbb{U}^+\ni w_k\searrow u^+$. Fix such a sequence $\{w_{k}\}_{k=1}^{\infty}$ and an arbitrary stochastic sub-solution $w_{-}$. Let $\tilde{\vp}(t,x)=\vp(t,x)+\iota|x-x_0|^{n_{0}}$.\footnote{Since we will fix $n_{0}$ and $\iota$ later, we still use the notation $\tilde\vp$ when without ambiguity despite the fact that the function depends on $n_{0}$ and $\iota$.} We can choose $n_{0}\geq 2$ such that for any $\iota>0$, 
\begin{equation}\label{eq:varphi grows faster}
\min_{0\leq t\leq T} (\tilde\vp(t,x)-w_{1}(t,x))\rightarrow\infty\;\; \text{as} \;\; |x|\rightarrow\infty.
\end{equation}
We can do this because $\vp(t, x)$ is bounded from below by $w_{-}$ (which has polynomial growth in $x$) and $w_{1}$ has polynomial growth in $x$. Since $(\mathcal{N}_{\eps,\eta})_{\eps\geq 0}$ is non-decreasing in $\eps$, we know 
\begin{equation*}
H_{*}(\Theta, \vp)=\liminf_{\begin{subarray}{c} \Theta^{'}\rightarrow\Theta, \psi \overset{\text{u.c.}}{\longrightarrow} \vp \\ \eta\searrow 0\end{subarray}} H_{0,\eta} (\Theta^{'}, \psi).
\end{equation*}
By \eqref{eq: HJB operators} and \eqref{eq: initial_local_contra_condition_at_t0x0_for u+}, we can find $\eps>0$, $\eta>0$ and $\iota>0$ such that for all $(t,x,y)$ satisfying $(t,x)\in B_{\eps}(t_0,x_0)$ and $|y-\tilde{\vp}(t,x)|\leq \eps$,
$
\mu_Y(t,x,y,u)-\mathcal{L}^u\tilde{\vp}(t,x)\geq 2\eta \text{ for some }u\in \mathcal{N}_{0, \eta}(t,x,y,D\tilde{\vp}(t,x),\tilde{\vp}).
$
Fix $\iota$. Note that $(t_0,x_0)$ is still a strict maximizer of $u^+-\tilde{\vp}$ over $\Di$. For $\eps$ sufficiently small, Assumption \ref{assump: regularity} implies that there exists a locally Lipschitz map $\hat{\nu}$ such that
\begin{gather}\label{eq: local_eqn1}
\hat{\nu}(t,x,y,D\tilde{\vp}(t,x))\in \mathcal{N}_{0,\eta}(t,x,y,D\tilde{\vp}(t,x),\tilde{\vp}) \text{ and} \\
\begin{array}{c} \label{eq: local_eqn2}
\mu_Y(t,x,y,\hat{\nu}(t,x,y,D\tilde{\vp}(t,x)))-\mathcal{L}^{\hat{\nu}(t,x,y,D\tilde{\vp}(t,x))}\tilde{\vp}(t,x)\geq \eta \\
\text{for all } (t,x,y)\in\Di\times \R\;\text{s.t.}\;(t,x)\in B_{\eps}(t_0,x_0) \text{ and }|y-\tilde{\vp}(t,x)|\leq \eps.
\end{array}
\end{gather}
In the arguments above, choose $\eps$ small enough such that $\text{cl}(B_{\eps}(t_0,x_0))\cap \DT=\emptyset$. Since \eqref{eq:varphi grows faster} holds,  there exists $ R_0>\eps$ such that $\tilde\vp>w_1+\eps\geq w_k+\eps$ on $\mathbb{O}:=\D\setminus[0,T]\times \text{cl}(B_{R_0}(x_0))\;\text{for all } k.$
On the compact set $\mathbb{T}:= [0,T]\times \text{cl}(B_{R_0}(x_0))\setminus B_{\eps/2}(t_0, x_0)$,  we know that $\tilde\varphi >u^+$ and the minimum of $\tilde\vp -u^+$ is attained since $u^+$ is USC. Therefore, $\tilde\vp >u^+ +2\alpha$ on $\mathbb{T}$ for some $\alpha >0$. By a Dini-type argument, for large enough $n$, we have $\tilde\vp >w_n+\alpha$ on $\mathbb{T}$ and $\tilde\vp>w_n-\eps$ on $\text{cl}(B_{\eps/2}(t_0, x_0))$. For simplicity, fix such an $n$ and set $w=w_n$. In short,
\begin{equation}\label{eq: comparion of w and varphi on differen sets_u+}
\tilde\vp>w+\eps\;\;\text{on}\;\; \mathbb{O},\; 
\tilde\vp >w+\alpha \;\;\text{on}\;\;\mathbb{T} \;\;\text{and}\;\; \tilde\vp>w-\eps\;\;\text{on}\;\; \text{cl}(B_{\eps/2}(t_0, x_0)).
\end{equation}
For $\kappa \in\;]0,\varepsilon\wedge\alpha[\;$, define
$$w^{\kappa}:=
\left \{ 
\begin{array}{ll}
(\tilde\varphi -\kappa)\wedge w\ \ {\textrm on}\ \  \text{cl}(B_{\eps}(t_0, x_0)),\\
w \ \ \textrm{outside}\ \ \text{cl}(B_{\eps}(t_0, x_0)).
\end{array}
\right.
$$
Observing that $w^{\kappa}(t_0,x_0)=\tilde\vp(t_0,x_0)-\kappa<u^{+}(t_0,x_0)$, we could obtain a contradiction if we could show that $w^{\kappa}\in\mathbb{U}^+$. Obviously, $w^{\kappa}$ is continuous, has polynomial growth in $x$ and $w^{\kappa}(T,x)\geq g(x)$ for all $x\in\mathbb{R}^d$. \\
Fix $(t,x,y)\in \Di\times \mathbb{R}$, $\nu\in \Uc^t$ and $\tau \in\mathcal{T}_t$.\footnote{Here we choose $(t,x)\in \Di$ since the case $(t,x)\in \DT$ is trivial.} Now our goal  is to construct an admissible control $\widetilde{\nu}$ such that $w^{\kappa}$ and the processes $(X,Y)$ controlled by $\nu\otimes_{\tau}\widetilde{\nu}$ satisfy the property in the definition of stochastic super-solutions.  \\
\indent Let $A=\{w^{\kappa}(\tau,X_{t,x}^{\nu}(\tau))= w(\tau,X_{t,x}^{\nu}(\tau))\}.$ On $A$, let $\widetilde{\nu}$ be $\widetilde{\nu}_1$, which is ``optimal" for $w$ starting at $\tau$. We get the existence of $\widetilde{\nu}_1$ since $w\in \mathbb{U}^+$. On $A^c$, by an argument similar to that in \cite{BayraktarLi} (see Step 1.1 of Theorem 3.1's proof), we can construct an admissible control $\nu_0\in\Uc^t$ such that
\begin{equation*}
\begin{array}{c}
\nu_0(s):=\hat{\nu}\left(s,X^{\nu\otimes_{\tau}\nu_0}_{t,x}(s),Y^{\nu\otimes_{\tau}\nu_0}_{t,x,y}(s),D\tilde\varphi(s, X^{\nu\otimes_{\tau}\nu_0}_{t,x}(s)\right) \;\; \text{for}\;\; \tau \leq s< \theta,\footnote{The control $\nu_0$ can be fixed outside the time interval $[\tau,\theta)$, since we are only interested in the restriction of it to $[\tau,\theta)$.}  \;\;\text{where}\;\; \theta=\theta_1\wedge\theta_2\;\; \text{and} \\
\theta_{1}:=\inf\left\{s \in [ \tau,T]: (s, X_{t,x}^{\nu \otimes_{\tau} \nu_0 }(s)) \notin B_{\eps/2}(t_0, x_0) \right\}\wedge T, \;\;
\theta_2:=\inf\left\{s \in [ \tau,T]: \left|Y_{t,x,y}^{\nu \otimes_{\tau} \nu_0}(s)-\tilde\varphi(s,X_{t,x}^{\nu \otimes_{\tau} \nu_0}(s))\right| \geq \varepsilon \right\}\wedge T.
\end{array}
\end{equation*}
In the construction of $\nu_0$, we take advantage of Assumption \ref{assump: regu_on_coeff} and the Lipschitz  continuity of $\hat{\nu}$ which guarantee the existence of $X^{\nu\otimes_{\tau}\nu_0}_{t,x}$ and  $Y^{\nu\otimes_{\tau}\nu_0}_{t,x,y}$. Since $X^{\nu\otimes_{\tau}\nu_0}_{t,x}$ and $Y_{t,x,y}^{\nu\otimes_{\tau}\nu_0}$ are c\`adl\`ag, it is easy to check that $\theta\in\mathcal{T}_{\tau}$. We also see that 
\begin{gather}
\label{eq: s,X,Y are note in a neighborhood for after theta for U+}
(\theta_1, X_{t,x}^{\nu \otimes_{\tau} \nu_0}(\theta_1))\notin B_{\eps/2}(t_0, x_0),\;\; \left|Y_{t,x,y}^{\nu \otimes_{\tau} \nu_0}(\theta_2)-\tilde\varphi(\theta_2,X_{t,x}^{\nu \otimes_{\tau} \nu_0}(\theta_2))\right|\geq \eps, \\
\label{eq: s,X,Y are in a neighborhood for before theta for U+}
(\theta_1, X_{t,x}^{\nu \otimes_{\tau} \nu_0}(\theta_1-))\in \text{cl}(B_{\eps/2}(t_0, x_0)),\;\; \left|Y_{t,x,y}^{\nu \otimes_{\tau} \nu_0}(\theta_2-)-\tilde\varphi(\theta_2,X_{t,x}^{\nu \otimes_{\tau} \nu_0}(\theta_2-))\right|\leq \eps.
\end{gather}
 Let $\widetilde{\nu}^{\theta}$ be the ``optimal" control for $w$ starting at $\theta$. We define $\tilde{\nu}$ on $A^c$ by $\nu_0\otimes_{\theta}\widetilde{\nu}^{\theta}$. In short,
\begin{equation*}
\widetilde{\nu}:=\left(\mathbbm{1}_{A}\widetilde{\nu}_1+\mathbbm{1}_{A^c}(\nu_0\mathbbm{1}_{[t,\theta[}+\mathbbm{1}_{[\theta,T]}\widetilde{\nu}^{\theta})\right)\mathbbm{1}_{[\tau,T]}.
\end{equation*}
It is not difficult to check that $\widetilde{\nu}\in\mathcal{U}^t$. To prove that the above construction works, we next show that
$Y(\rho )\geq w^{\kappa}(\rho, X(\rho ))\;\; \mathbb{P}-\text{a.s.}$ on $\{Y(\tau )\geq w^{\kappa}(\tau, X(\tau))\},$
where  $X:= X_{t,x}^{\nu\otimes_{\tau}\tilde{\nu}}$ and $Y:=Y_{t,x,y}^{\nu\otimes_{\tau}\tilde{\nu}}$.  Corresponding to the construction of $\widetilde{\nu}$ on $A$ and $A^{c}$, we consider the following two cases: \\
\noindent\textbf{(i) On the set $A\cap\{Y(\tau)\geq w^{\kappa}(\tau,X(\tau))\}$}. We have 
$
Y(\tau)\geq w(\tau,X(\tau)).
$
From the definition of $\nu$ on $A$ and the fact that $w\in\bU^{+}$, we know  
$$
Y(\rho)=Y_{t,x,y}^{\nu\otimes_{\tau}\widetilde{\nu}_1}(\rho)\geq w(\rho,X_{t,x}^{\nu\otimes_{\tau}\widetilde{\nu}_1}(\rho))\geq w^{\kappa}(\rho,X(\rho)) \;\; \mathbb{P}-\text{a.s}\;\; \text{on}\;\;A\cap\{Y(\tau)\geq w^{\kappa}(\tau,X(\tau))\}.
$$
\textbf{(ii) On the set $A^c \cap\{Y(\tau)\geq w^{\kappa}(\tau,X(\tau))\}$}. Letting $\Gamma(s):=Y(s)-\tilde\vp(s,X(s))$,  we use It\^{o}'s formula and the definition of $\nu_0$ to obtain
\begin{equation*}
\Gamma(\cdot \wedge \theta)=\Gamma(\tau) +\int_{\tau}^{\cdot\wedge\theta}\int_{E}\overline{J}^{\nu_0(s),e}\left(s,Z(s-),\tilde\vp\right) ^{\top}\lambda(ds,de)+ \int_{\tau}^{\cdot\wedge\theta}\left(\mu_Y(s, Z(s), \nu_0(s))-\mathscr{L}^{\nu_0(s)}\tilde{\vp}(s, X(s))\right)ds
\end{equation*}
on $A\cap\{Y(\tau)\geq w^{\kappa}(\tau,X(\tau))\}.$ Therefore, by \eqref{eq: local_eqn1}, \eqref{eq: local_eqn2}, \eqref{eq: s,X,Y are in a neighborhood for before theta for U+} and the definition of $\theta$,  we know that $\Gamma(\cdot \wedge \theta)$ is non-decreasing on $[\tau,T]$. This implies that 
\begin{equation}\label{eq: increasing_process_tau_theta}
Y(\theta)- \tilde \vp(\theta,X(\theta))+\kappa\geq Y(\tau)-  \tilde\vp(\tau,X(\tau))+\kappa\geq 0\;\;\text{on}\;\;A^c \cap\{Y(\tau)\geq w^{\kappa}(\tau,X(\tau))\}.
\end{equation}
Since $(\theta_1, X(\theta_1))\notin B_{\eps/2}(t_0, x_0)$, we know
\begin{equation}\label{eq:theta1_less_theta2}
0\leq Y(\theta_1)-  \tilde\vp(\theta_1,X(\theta_1))+\kappa\leq Y(\theta_1)- w(\theta_1,X(\theta_1)) \;\; \text{on} \;\; \{\theta_1\leq\theta_2\}\cap A^c \cap\{Y(\tau)\geq w^{\kappa}(\tau,X(\tau))\}
\end{equation}
from \eqref{eq: comparion of w and varphi on differen sets_u+}. On the other hand, it holds that 
$Y(\theta_2)-\tilde\vp(\theta_2,X(\theta_2)) \geq \eps $ on $ \{\theta_1>\theta_2\}\cap A^c \cap\{Y(\tau)\geq w^{\kappa}(\tau,X(\tau))\}$  due to \eqref{eq: s,X,Y are note in a neighborhood for after theta for U+} and \eqref{eq: increasing_process_tau_theta}. Therefore, since $\tilde\varphi>w-\eps$ on $\text{cl}(B_{\eps/2}(t_0, x_0))$ and \eqref{eq: s,X,Y are in a neighborhood for before theta for U+} holds,
\begin{equation}\label{eq:theta1_greater_theta2}
 Y(\theta_2)- w(\theta_2,X(\theta_2))\geq\eps+\tilde\varphi(\theta_2,X(\theta_2))-w(\theta_2,X(\theta_2))>0 \;\; \text{on}\;\; \{\theta_1>\theta_2\}\cap A^c \cap\{Y(\tau)\geq w^{\kappa}(\tau,X(\tau))\}.
 \end{equation}
 Combining \eqref{eq:theta1_less_theta2} and \eqref{eq:theta1_greater_theta2}, we obtain
$Y(\theta)- w(\theta,X(\theta))\geq 0$ on $A^c \cap\{Y(\tau)\geq w^{\kappa}(\tau,X(\tau))\}.$
Therefore, from the definition of $\widetilde{\nu}^{\theta}$, 
\begin{equation}\label{eq:aits}
 Y(\rho \vee \theta)- w^{\kappa}(\rho \vee \theta,X(\rho \vee \theta))  \geq  Y(\rho \vee \theta)- w(\rho \vee \theta,X(\rho \vee \theta))  \geq 0\;\;\text{on}\;\; A^c \cap\{Y(\tau)\geq w^{\kappa}(\tau,X(\tau))\}.
\end{equation}
Also, the monotonicity of $ \Gamma(\cdot \wedge \theta)$ implies that $Y(\rho \wedge \theta)-\tilde\vp(\rho \wedge \theta,X(\rho\wedge \theta))+\kappa\geq 0$ on $A^c \cap\{Y(\tau)\geq w^{\kappa}(\tau,X(\tau))\}.$
This means that 
\begin{equation}\label{eq:eans}
\mathbbm{1}_{\{\rho<\theta\}}\left(Y(\rho)-w^{\kappa}(\rho,X(\rho))\right)\geq 0\;\;\text{on}\;\;  A^c \cap\{Y(\tau)\geq w^{\kappa}(\tau,X(\tau))\}.
\end{equation}
From \eqref{eq:aits} and \eqref{eq:eans}, we get $Y(\rho)-w^{\kappa}(\rho,X(\rho))\geq 0 \;\; \text{on} \,\,  A^c \cap\{Y(\tau)\geq w^{\kappa}(\tau,X(\tau))\}.$

{\bf \noindent Step 2 ($u^-$ is a viscosity super-solution).}
Let $(t_0,x_0)\in\Di$ satisfy $0=(u^--\varphi)(t_0,x_0)= \min_{\Di}(u^--\varphi)$
for some $\vp\in C^{1,2}(\D)$. For the sake of contradiction, assume that
\begin{equation}\label{eq: initial_local_contra_condition_at_t0x0_for u-}
-2\eta:=-\partial_t\vp(t_0,x_0)+H^*\vp(t_0,x_0)<0.
\end{equation}
Let $\{w_{k}\}_{k=1}^{\infty}$ be a sequence in $\bU^{-}$ such that $ w_k\nearrow u^-$. Let $\tilde{\vp}(t,x):=\vp(t,x)-\iota|x-x_0|^{n_{0}}$, where we choose $n_{0}\geq 2$ such that for all $\iota>0$,
\begin{equation}\label{eq: varphi decays faster}
 \max_{0\leq t\leq T}(\tilde\vp(t,x)-w_{1}(t,x))\rightarrow-\infty\;\; \text{and}\;\;  \max_{0\leq t\leq T} \tilde\vp(t,x)\rightarrow-\infty\;\;\text{as}\;\; |x|\rightarrow \infty. \footnote{The existence of $n_{0}$ follows as in Step1.}
\end{equation}
By \eqref{eq: initial_local_contra_condition_at_t0x0_for u-}, the upper semi-continuity of $H^*$ and the fact that $\tilde{\vp}\overset{\text{u.c.}}{\longrightarrow}\vp$ as $\iota\rightarrow 0$, we can find $\eps>0$, $\eta>0$ and $\iota>0$ such that
\begin{equation}\label{eq: local contra u^-}
\begin{array}{c}
\mu_Y(t,x,y,u)-\mathscr{L}^{u}\tilde{\vp}(t,x)\leq-\eta\;\text{for all }u\in\mathcal{N}_{\eps, -\eta}(t,x,y,D\tilde{\vp}(t,x),\tilde{\vp}) \\
\text{and}\;\;(t,x,y)\in\Di\times\R\;\text{s.t.}\;(t,x)\in B_{\eps}(t_0,x_0)\;\text{and}\;|y-\tilde{\vp}(t,x)|\leq \eps.
\end{array}
\end{equation}
Fix $\iota$. Note that $(t_0,x_0)$ is still a strict minimizer of $u^--\tilde{\vp}$.  Since \eqref{eq: varphi decays faster} holds, there exists $R_0>\eps$ such that $$\tilde\vp<w_1-\eps\leq w_k-\eps\;\;\text{on}\;\; \mathbb{O}:=\D\setminus[0,T]\times \text{cl}(B_{R_0}(x_0)).$$
On the compact set $\mathbb{T}:= [0,T]\times \text{cl}(B_{R_0}(x_0))\setminus B_{\eps/2}(t_0, x_0)$,  we know that $\tilde\vp <u^-$ and the maximum of $\tilde\vp-u^-$ is attained since $u^-$ is LSC. Therefore, $\tilde\vp <u^- -2\alpha$ on $\mathbb{T}$ for some $\alpha >0$. By a Dini-type argument, for large enough $n$, we have $\tilde\varphi <w_n-\alpha$ on $\mathbb{T}$ and $\tilde\varphi<w_n+\varepsilon$ on $\text{cl}(B_{\eps/2}(t_0, x_0))$. For simplicity, fix such an $n$ and set $w=w_n$. In short,
\begin{equation}\label{eq: comparion of w and varphi on differen sets_u-}
\tilde\vp<w-\eps\;\;\text{on}\;\; \mathbb{O},\; 
\tilde\varphi <w-\alpha \;\;\text{on}\;\;\mathbb{T} \;\;\text{and}\;\; \tilde\varphi<w+\varepsilon\;\;\text{on}\;\; \text{cl}(B_{\eps/2}(t_0, x_0)).
\end{equation}
For $\kappa\in\;]0,\alpha\wedge \varepsilon[\;$, define
$$w^{\kappa}:=
\left \{ 
\begin{array}{ll}
(\tilde\varphi +\kappa)\vee w\ \ {\textrm on}\ \  \text{cl}(B_{\eps}(t_0, x_0)),\\
w \ \ \textrm{outside}\ \  \text{cl}(B_{\eps}(t_0, x_0)).
\end{array}
\right.
$$
Noticing that $w^{\kappa}(t_0,x_0)\geq\tilde\varphi(t_0,x_0)+\kappa>u^-(t_0,x_0)$, we will obtain a contradiction if we show that $w^{\kappa}\in \mathbb{U}^-$. Obviously, $w^{\kappa}$ is continuous, has polynomial growth in $x$ and $w^{\kappa}(T,x)\leq g(x)$ for all $x\in\mathbb{R}^d$. Fix $(t,x,y)\in \Di\times \mathbb{R}$, $\nu\in \Uc^t$ and $\tau \in\mathcal{T}_t$. Our goal is to show that
$$\mathbb{P}(Y(\rho )< w^{\kappa}(\rho, X(\rho ))|B)>0$$
for all $\rho\in\mathcal{T}_{\tau}$ and $B\subset \{Y(\tau)<w^{\kappa}(\tau,X(\tau))\}$ satisfying $B\in\mathcal{F}_\tau^t$ and $\P(B)>0$, where $X:= X_{t,x}^{\nu}$ and $Y:=Y_{t,x,y}^{\nu}.$ Let $A=\{w^{\kappa}(\tau,X(\tau))=w(\tau,X(\tau))\}$ and set
\begin{equation*}
\begin{array}{c}
E=\{Y(\tau)< w^{\kappa}(\tau,X(\tau))\}, \;\;
E_0=E\cap A, \;\; E_1=E\cap A^c,  \\
G=\{Y(\rho)< w^{\kappa}(\rho,X(\rho)\}, \;\;
G_0=\{Y(\rho)< w(\rho,X(\rho)\}.
\end{array}
\end{equation*}
Then $E=E_0\cup E_1,\; E_0\cap E_1=\emptyset \;\;\text{and}\;\; G_0\subset G. $
To prove that $w^{\kappa}\in\bU^{-}$, it suffices to show that $\mathbb{P}(G\cap B)>0$. 
As in \cite{BayraktarLi} and \cite{SONER_TOUZI_STG}, we will show $\P(B\cap E_{0})>0\implies \P(G\cap B\cap E_{0})>0$ and $\P(B\cap E_{1})>0\implies \P(G\cap B\cap E_{1})>0$. This, together with the facts $\P(B)=\P(B\cap E_{0})+\P(B\cap E_{1})>0$ and $\P(G\cap B)=\P(G\cap B\cap E_{0})+\P(G\cap B\cap E_{1})$, implies that $\P(G\cap B)>0$.\\
\textbf{(i)}{\bf Assume that $\mathbb{P}(B\cap E_0)>0$}.  Since $B\cap E_0 \subset \{Y(\tau)<w(\tau,X(\tau))\}$ and $B\cap E_0 \in \mathcal{F}_\tau^t$, $\mathbb{P}(G_0|B\cap E_0)>0$ from the definition of $\bU^{-}$. This further implies that $\mathbb{P}(G\cap B\cap E_0) \geq \mathbb{P}(G_0\cap B\cap E_0)>0$. \\
\textbf{(ii)}{\bf Assume that $\mathbb{P}(B\cap E_1)>0$}. Let $\theta=\theta_1\wedge\theta_2$, where
\begin{equation*}
\theta_{1}:=\inf\left\{s \in [ \tau,T]: (s, X(s)) \notin B_{\eps/2}(t_0, x_0) \right\}\wedge T, \;\; \theta_2:=\inf\left\{s \in [ \tau,T]: \left|Y(s)-\tilde\varphi(s,X(s))\right| \geq \varepsilon \right\} \wedge T.
\end{equation*}
Since $X$ and $Y$ are c\`adl\`ag processes, we know that $\theta\in\mathcal{T}_{\tau}$. The following also hold:  
\begin{gather}
\label{eq: theta-, X(theta-), Y(theta-) not are in neighborhood}
(\theta_1, X(\theta_1))\notin B_{\eps/2}(t_0, x_0),\;\;
\left|Y(\theta_2)-\tilde\varphi(\theta_2,X(\theta_2))\right|\geq \eps, \\
\label{eq: theta-, X(theta-), Y(theta-) are in neighborhood}
(\theta_1, X(\theta_1-))\in \text{cl}(B_{\eps/2}(t_0, x_0)),\;\;
\left|Y(\theta_2-)-\tilde\varphi(\theta_2,X(\theta_2-))\right|\leq \eps.
\end{gather}
Let
\begin{equation*}
\begin{gathered}
c^e_i(s) = J_i^{u,e}(s,X(s-), Y(s-), \tilde{\vp}), \;\;
d_i(s)=\int_E c^e_i(s) m_i(de), \;\; d(s)=\sum_{i=1}^{I}d_i(s), \\
a(s) = \mu_Y(s, X(s), Y(s), \nu(s))-\mathscr{L}^{\nu(s)}\tilde{\vp}(s, X(s)),\;\;
\pi(s)=N^{\nu(s)}(s, X(s), Y(s), D\tilde{\vp}(s, X(s))), \\
A_0 = \left\{s\in[\tau, \theta]: |\pi(s)|\leq \eps\right\}, \;\; A_{3,i}=\left\{(s,e)\in[\tau,\theta]\times E: c_i^e(s)\leq -\eta/2 \right\},  \\
A_1 = \left\{s\in[\tau, \theta]: c_i^e(s)\geq -\eta\text{ for } \hat{m}-a.s.\;e\in E  \text{ for all } i=1,\cdots, I \right\},\;\;  A_2 = (A_1)^c.
\end{gathered}
\end{equation*}
We then set
\begin{equation*}
L(\cdot):=\mathcal{E}\left(\int_t^{\cdot\wedge\theta}\int_E \sum\delta_i^e(s)\tilde\lambda_i(ds,de)+\int_{t}^{\cdot\wedge\theta}\alpha^{\top}(s) dW_s\right),
\end{equation*} 
where $\mathcal{E}(\cdot)$ denotes the Dol\'eans-Dade exponential and
\begin{equation*}
\begin{array}{c}
x^{+}:=\max\{0,x\},\;\;x^{-}:=\max\{0,-x\},\;\; \alpha(s):=-\frac{a(s)+d(s)}{|\pi(s)|^2}\pi(s)\mathbbm{1}_{A_0^c}(s), \;\;
M_i(s):=\int_E\mathbbm{1}_{A_{3,i}}(s,e)m_i(de), \\ 
K_i(s,e):=\left \{ 
\begin{array}{ll}
\frac{\mathbbm{1}_{A_{3,i}}(s,e)}{M_i(s)} \;\;\text{if}\;\; M_i(s)=0 \\
0 \quad\quad\quad\quad \text{otherwise}
\end{array}
\right.,
\;\;\delta_i^e(s):=\left(\frac{\eta}{2(1+|d(s)|)}-1+\mathbbm{1}_{A_2}(s)\cdot\frac{2a(s)^++\eta}{\eta }\cdot K_i(s,e)\right)\mathbbm{1}_{A_0}(s).
\end{array}
\end{equation*} 
If $s\in A_2$, then it follows from Assumption \ref{assump: lambda intensity kernel} and definitions of $A_2$ and $A_{3,i}$ that
\begin{equation}\label{eq: at least A_3i is not empty}
M_{i_0}(s) >0 \text{ for some }i_0\in\{1,2,\cdots, I\}. 
\end{equation}
Obviously, $L$ is a nonnegative local martingale on $[t,T]$. Therefore, it is a super-martingale. Let $\Gamma(s):=Y(s)-\tilde{\vp}(s,X(s))-\kappa$. Applying It\^o's formula, we get
\begin{eqnarray*}
 \Gamma(\cdot\wedge\theta)L(\cdot\wedge\theta)= \Gamma(\tau)L(\tau)+ \int_{\tau}^{\cdot\wedge\theta} L(s)\left\{\left(a(s)+d(s)\right)\mathbbm{1}_{A_0}(s)+\int_E\sum c^e_i(s)\delta^e_i(s)m_i(de)\right\}ds \\ 
 \int_{\tau}^{\cdot\wedge\theta}\int_E\sum L(s)\left\{c^e_i(s)+\Gamma(s)\delta^e_i(s)+c^e_i(s)\delta^e_i(s)\right\}\tilde\lambda(ds,de)
+\int_{\tau}^{\cdot\wedge\theta}L(s)\left(\pi(s)+ \Gamma(s)\alpha(s)\right)^{\top}dW_s.
\end{eqnarray*} 
By the definition of $\delta_i^e$ and the fact that $\mathbbm{1}_{A_1}+\mathbbm{1}_{A_2}=1$ on $[\tau, \theta]$, the first integral in the equation above is 
\begin{eqnarray*}
\begin{split}
 & \displaystyle\int_{\tau}^{\cdot\wedge\theta} L(s) \left\{ \left(a(s)+\frac{\eta d(s)}{2(|d(s)|+1)}\right)\mathbbm{1}_{A_0\cap A_1}(s)+\mathbbm{1}_{A_0\cap A_2}(s) \right. \\
& \times  \left.  \left(a(s)+\frac{\eta d(s)}{2(|d(s)|+1)}+ \frac{2a(s)^++\eta}{\eta }\int_E\sum c^e_i(s)K_i(s,e)m_i(de)\right) \vphantom{\frac{\eta d(s)}{2(|d(s)|+1)}}\right\} ds.
\end{split}
\end{eqnarray*} 
By \eqref{eq: local contra u^-}, $a(s)\leq -\eta$ on $A_0\cap A_1$. Then,
\begin{equation}\label{eq: drift1<0 for u-}
\left(a(s)+\frac{\eta d(s)}{2(|d(s)|+1)}\right)\mathbbm{1}_{A_0\cap A_1}(s) \leq \left(-\eta+\frac{\eta}{2}\right)\mathbbm{1}_{A_0\cap A_1}(s)\leq 0.
\end{equation}
By the definition of $A_{3,i}$ and \eqref{eq: at least A_3i is not empty}, it holds that
\begin{equation}\label{eq: drift2<0 for u-}
\begin{split}
& \mathbbm{1}_{A_0\cap A_2}(s)\left(a(s)+\frac{\eta d(s)}{2(|d(s)|+1)}+ \frac{2a(s)^++\eta}{\eta }\int_E\sum c^e_i(s)K_i(s,e)m_i(de)\right)\\
\leq & \mathbbm{1}_{A_0\cap A_2}(s)\left(a(s)+\frac{\eta}{2}-\frac{2a(s)^++\eta}{\eta }\cdot\frac{\eta}{2}\right)=- \mathbbm{1}_{A_0\cap A_2}(s)a(s)^-.
\end{split}
\end{equation}
Therefore, \eqref{eq: drift1<0 for u-} and \eqref{eq: drift2<0 for u-} imply that $\Gamma L$ is a local super-martingale on $[\tau,\theta]$. 
Note that
\begin{equation*}
\Gamma(\theta)-\Gamma(\theta-)=\int_{E} \overline{J}^{\nu(\theta),e}\left(\theta,X(\theta-),Y(\theta-),\tilde\vp\right) ^{\top}\lambda(\{\theta\},de).
\end{equation*}
Since $\tilde\vp\in C(\D)$ and \eqref{eq: varphi decays faster} holds, $\tilde\vp$ is locally bounded and globally bounded from above. This, together with \eqref{eq: theta-, X(theta-), Y(theta-) are in neighborhood} and the admissibility condition \eqref{eq:admissibility}, 
implies that $
\Gamma(\theta)-\Gamma(\theta-)\geq -K$ almost surely for some $K>0$ ($K$ may depend on $(t_{0}, x_{0}), \eps$, $\nu$ and $\tilde\vp$). Since
$
\Gamma(s)= Y(s)-\tilde{\vp}(s,X(s))-\kappa\geq -(\eps+\kappa) \text{ on }[\tau,\theta[\;
$, $\Gamma L$ is bounded from below by a sub-martingale $-(\eps+\kappa+K)L$ on $[\tau, \theta]$. This further implies that $\Gamma L$ is a super-martingale by Fatou's Lemma. Since $\Gamma(\tau)L(\tau)<0$ on $B\cap E_1$, the super-martingale property implies that there exists $F \subset B \cap E_1$ such that $F\in\mathcal{F}^t_{\tau}$ and $\Gamma(\theta\wedge \rho)L(\theta\wedge\rho)<0$ on $F$. The non-negativity of $L$ then yields $\Gamma(\theta\wedge\rho)<0$.  Therefore, 
\begin{gather}
Y(\theta_1)<\tilde\varphi(\theta_1,X(\theta_1))+\kappa \text{ on }F\cap \{\theta_1\leq\theta_2, \theta<\rho \},\;\;
Y(\theta_2)<\tilde\varphi(\theta_2,X(\theta_2))+\kappa\text{ on }F\cap \{\theta_1>\theta_2, \theta<\rho\}\text{ and} \nonumber \\
\label{eqn_thirdpart_of_delta}
Y(\rho)-(\tilde\varphi(\rho ,X(\rho))+\kappa)<0 \text{ on } F\cap \{\theta\geq \rho\}.
\end{gather}
Since $(\theta_1,X(\theta_1))\notin B_{\eps/2}(t_0,x_0)$, it follows from the first two inequalities in \eqref{eq: comparion of w and varphi on differen sets_u-} that
\begin{equation}\label{eq: case theta1<theta2,u-}
 Y(\theta_1)<\tilde\vp(\theta_1,X(\theta_1))+\kappa<w(\theta_1,X(\theta_1))\;\;\text{on}\;\; F\cap \{\theta_1\leq\theta_2, \theta<\rho \}. 
\end{equation}
On the other hand, since $Y(\theta_2)<\tilde\varphi(\theta_2,X(\theta_2))+\kappa$ on $ F\cap \{\theta_1>\theta_2, \theta<\rho \}$ and \eqref{eq: theta-, X(theta-), Y(theta-) not are in neighborhood} holds, 
$Y(\theta_2)-\tilde\vp(\theta_2,X(\theta_2))\leq -\eps$ on $F\cap \{\theta_1>\theta_2, \theta<\rho \}$. Observing that $(\theta_2,X(\theta_2))\in B_{\eps/2}(t_0,x_0)$ on $\{\theta_{1}>\theta_{2}\}$, we get from the last inequality of \eqref{eq: comparion of w and varphi on differen sets_u-} that
\begin{equation}\label{eq: case theta1>theta2, u-}
Y(\theta_2)-w(\theta_2,X(\theta_2))<\tilde\vp(\theta_2,X(\theta_2))-\eps-w(\theta_2,X(\theta_2))< 0 \;\;\text{on}\;\; F\cap \{\theta_1>\theta_2, \theta<\rho\}.
\end{equation}
From \eqref{eq: case theta1<theta2,u-} and \eqref{eq: case theta1>theta2, u-}, we get that $Y(\theta)<w(\theta,X(\theta))$ on $F\cap\{\theta< \rho\}.$ Therefore, from the definition of $\bU^{-}$, 
\begin{equation} \label{eq: conlcusion_on_theta_less_than_rho}
\mathbb{P}(G_0|F\cap\{\theta< \rho\})>0 \;\; \text{ if} \;\; \mathbb{P}(F\cap\{\theta< \rho\})>0.
\end{equation} 
From \eqref{eqn_thirdpart_of_delta}, it holds that
\begin{equation} \label{eq: conlcusion_on_theta_larger_than_rho}
\mathbb{P}(G|F\cap\{\theta\geq \rho\})>0 \;\; \text{ if} \;\; \mathbb{P}(F\cap\{\theta\geq \rho\})>0.
\end{equation}
Since $G_0\subset G$, \eqref{eq: conlcusion_on_theta_less_than_rho} and \eqref{eq: conlcusion_on_theta_larger_than_rho} imply that
$\mathbb{P}(G\cap F)>0.$
Therefore, $\mathbb{P}(G\cap B\cap E_1)>0$. \qed
\end{singlespace}

\begin{singlespace}
\small 
\noindent \textbf{Acknowledgments} \;\;We would like to thank Bruno Bouchard who encouraged us to write this paper and for his constructive comments on its first version. We also thank the referees and the anonymous associate editor for their helpful comments, which helped us to improve our paper.  This research is supported in part by the National Science Foundation.
\end{singlespace}

\bibliographystyle{spmpsci_unsrt}
\bibliography{mybib_T}

\begin{thebibliography}{10}
\providecommand{\url}[1]{{#1}}
\providecommand{\urlprefix}{URL }
\expandafter\ifx\csname urlstyle\endcsname\relax
  \providecommand{\doi}[1]{DOI~\discretionary{}{}{}#1}\else
  \providecommand{\doi}{DOI~\discretionary{}{}{}\begingroup
  \urlstyle{rm}\Url}\fi

\bibitem{Soner_Touzi_Superreplication}
Soner, H.M., Touzi, N.: Superreplication under gamma constraints.
\newblock SIAM J. Control Optim. \textbf{39}(1), 73--96 (2000)

\bibitem{DP_FOR_STP_AND_G_FLOW}
Soner, H.M., Touzi, N.: Dynamic programming for stochastic target problems and
  geometric flows.
\newblock J. Eur. Math. Soc. (JEMS) \textbf{4}(3), 201--236 (2002)

\bibitem{SONER_TOUZI_STG}
Soner, H.M., Touzi, N.: Stochastic target problems, dynamic programming, and
  viscosity solutions.
\newblock SIAM J. Control Optim. \textbf{41}(2), 404--424 (2002)

\bibitem{Bruno_jump_diffusion}
Bouchard, B.: Stochastic targets with mixed diffusion processes and viscosity
  solutions.
\newblock Stochastic Process. Appl. \textbf{101}(2), 273--302 (2002)

\bibitem{Moreau}
Moreau, L.: Stochastic target problems with controlled loss in jump diffusion
  models.
\newblock SIAM Journal on Control and Optimization \textbf{49}(6), 2577--2607
  (2011)

\bibitem{Bouchard_Elie_Touzi_ControlledLoss}
Bouchard, B., Elie, R., Touzi, N.: Stochastic target problems with controlled
  loss.
\newblock SIAM J. Control Optim. \textbf{48}(5), 3123--3150 (2009/10)

\bibitem{Bayraktar_and_Sirbu_SP_LinearCase}
Bayraktar, E., S{\^{\i}}rbu, M.: Stochastic {P}erron's method and verification
  without smoothness using viscosity comparison: the linear case.
\newblock Proc. Amer. Math. Soc. \textbf{140}(10), 3645--3654 (2012)

\bibitem{Bayraktar_and_Sirbu_SP_DynkinGames}
Bayraktar, E., S{\^{\i}}rbu, M.: Stochastic {P}erron's method and verification
  without smoothness using viscosity comparison: obstacle problems and {D}ynkin
  games.
\newblock Proc. Amer. Math. Soc. \textbf{142}(4), 1399--1412 (2014)

\bibitem{Bayraktar_and_Sirbu_SP_HJBEqn}
Bayraktar, E., S{\^{\i}}rbu, M.: Stochastic {P}erron's method for
  {H}amilton-{J}acobi-{B}ellman equations.
\newblock SIAM J. Control Optim. \textbf{51}(6), 4274--4294 (2013)

\bibitem{MR3488161}
Claisse, J., Talay, D., Tan, X.: A {P}seudo-{M}arkov {P}roperty for
  {C}ontrolled {D}iffusion {P}rocesses.
\newblock SIAM J. Control Optim. \textbf{54}(2), 1017--1029 (2016)

\bibitem{Bouchard_Equivalence}
Bouchard, B., Dang, N.M.: Optimal control versus stochastic target problems: an
  equivalence result.
\newblock Systems Control Lett. \textbf{61}(2), 343--346 (2012)

\bibitem{BayraktarLi}
Bayraktar, E., Li, J.: Stochastic {P}erron for stochastic target games.
\newblock Ann. Appl. Probab. \textbf{26}(2), 1082--1110 (2016)

\bibitem{Barle-IntegroPDE}
Barles, G., Imbert, C.: Second-order elliptic integro-differential equations:
  viscosity solutions' theory revisited.
\newblock Ann. Inst. H. Poincar\'e Anal. Non Lin\'eaire \textbf{25}(3),
  567--585 (2008)

\bibitem{ShreveKaratzas}
Karatzas, I., Shreve, S.E.: Brownian motion and stochastic calculus,
  \emph{Graduate Texts in Mathematics}, vol. 113, second edn.
\newblock Springer-Verlag, New York (1991)

\bibitem{Bouchard_Nutz_TargetGames}
Bouchard, B., Nutz, M.: Stochastic target games and dynamic programming via
  regularized viscosity solutions.
\newblock Mathematics of Operations Research \textbf{41}(1), 109--124 (2016)

\end{thebibliography}
\end{document}